\documentclass[ijoo,nonblindrev]{informs-anon}

\OneAndAHalfSpacedXI

\usepackage{adjustbox}
\usepackage{mathtools}
\usepackage{xcolor}         
\usepackage{enumitem}
\usepackage[textsize=tiny, color=yellow]{todonotes}

\usepackage{tikz}
\usepackage{amsmath}
\usepackage{listings}
\usepackage{float}
\usepackage{subcaption}
\usepackage{multirow}
\usepackage{booktabs} 
\usepackage[ruled,vlined,linesnumbered]{algorithm2e}
\SetKwData{Left}{left}\SetKwData{This}{this}\SetKwData{Up}{up}                
\SetKwInput{KwInput}{Input}                
\SetKwInput{KwOutput}{Output}              
\SetKwIF{If}{ElseIf}{Else}{if}{}{else if}{else}{end if}%
\SetKwFor{While}{while}{}{end while}%
\SetKwRepeat{Do}{do}{while}
\SetKw{KwGoTo}{go to}
\usepackage{ifthen}  
\usepackage[hidelinks]{hyperref}
\usetikzlibrary{shapes,arrows}
\graphicspath{{figures/}} 

\makeatletter 
\newcommand{\owntag}[2][\relax]{
	\ifx#1\relax\relax\def\owntag@name{#2}\else\def\owntag@name{#1}\fi
	\refstepcounter{equation}\tag{\theequation, #2}%
	\expandafter\ltx@label\expandafter{eq:\owntag@name}%
	\edef\@currentlabel{\theequation, #2}\expandafter\ltx@label\expandafter{Eq:\owntag@name}%
	\def\@currentlabel{#2}\expandafter\ltx@label\expandafter{tag:\owntag@name}%
}
\makeatother 
\newcommand{\iter}{\operatorname{it}}
\newcommand{\ssclrmcns}{Alt-PZF-EnSC+gLRMC}
\newcommand{\ssclrmc}{{\ssclrmcns} }
\newcommand{\ourmethodns}{MISS-DSG}
\newcommand{\ourmethod}{{\ourmethodns} }

\newcommand{\ed}{ED}
\newcommand{\numsub}{K}
\newcommand{\cind}{i}
\newcommand{\ambientdimension}{d}
\newcommand{\numinitbasis}{b}
\newcommand{\numvectors}{n}
\newcommand{\numcols}{T}
\newcommand{\column}{t}
\newcommand{\datavector}{j}
\newcommand{\regularizer}{\lambda}
\newcommand{\missingpercentage}{f}
\newcommand{\distance}[1][]{
	\ifthenelse{ \equal{#1}{} }
	{\Delta_{\datavector\column}} 
	{\Delta_{#1}} 
}
\newcommand{\cost}[1][]{
	\ifthenelse{ \equal{#1}{} }
	{c_{\datavector\column}} 
	{c_{#1}} 
}

\newcommand{\assignvar}[1][]{
	\ifthenelse{ \equal{#1}{} }
	{x_{\datavector\column}} 
	{x_{#1}} 
}

\newcommand{\selectvar}[1][]{
	\ifthenelse{ \equal{#1}{} }
	{z_{\column}} 
	{z_{#1}} 
}

\newcommand{\costvar}[1][]{
	\ifthenelse{ \equal{#1}{} }
	{w_{\datavector}} 
	{w_{#1}} 
}

\newcommand{\observedset}[1][]{
	\ifthenelse{ \equal{#1}{} }
	{\Omega} 
	{\Omega_#1} 
}

\newcommand{\lowrank}[1][]{
	\ifthenelse{ \equal{#1}{} }
	{r} 
	{r_{#1}} 
}

\newcommand{\basis}[1][]{
	\ifthenelse{ \equal{#1}{} }
	{U} 
	{U_{#1}} 
}

\newcommand{\datasubset}[1][]{
	\ifthenelse{ \equal{#1}{} }
	{X(S)} 
	{X({#1})} 
}
\newcommand{\projectionOperator}[1][]{
	\ifthenelse{ \equal{#1}{} }
	{P_{\basis_{\observedset,\datavector}}} 
	{P_{#1}} 
}

\usepackage{natbib}
 \bibpunct[, ]{(}{)}{,}{a}{}{,}%
 \def\bibfont{\small}%
\TheoremsNumberedThrough     
\ECRepeatTheorems

\EquationsNumberedThrough    

\MANUSCRIPTNO{}

\graphicspath{{figures/}} 

\begin{document}



 \RUNAUTHOR{Soni et al.} 

\RUNTITLE{Subspace Clustering with Integer Programming}

\TITLE{
An Integer Programming Approach To Subspace Clustering With Missing Data}

\ARTICLEAUTHORS{%
\AUTHOR{Akhilesh Soni, Jeff Linderoth, James Luedtke}
\AFF{Department of Industrial and Systems Engineering, University of Wisconsin-Madison, \EMAIL{akhileshsoni95@gmail.com},
\EMAIL{linderoth@wisc.edu},\EMAIL{jim.luedtke@wisc.edu}} 
\AUTHOR{Daniel Pimentel-Alarc{\'o}n}
\AFF{Department of Biostatistics and Medical Informatics, University of Wisconsin-Madison , \EMAIL{pimentelalar@wisc.edu}} 
} 

\ABSTRACT{%
In the \emph{Subspace Clustering with Missing Data} (SCMD) problem, we are given a collection of $\numvectors$ partially
observed $\ambientdimension$-dimensional vectors.
The data points are assumed to be concentrated near a union of low-dimensional subspaces.
The goal of SCMD is to cluster  the vectors  according to their subspace membership and recover the underlying basis,
which can then be used to infer their missing entries.
State-of-the-art algorithms for SCMD 
can fail on instances with a high proportion of missing data, full-rank data, or if the underlying subspaces are similar
to each other.
We propose a novel integer programming approach for SCMD.
The approach is based on dynamically determining a set of candidate subspaces and optimally assigning points to selected subspaces.  The problem structure is identical to the classical facility-location problem, with subspaces playing the role of facilities and data points that of customers.  We propose a column-generation approach for identifying candidate subspaces combined with a Benders decomposition approach for solving the linear programming relaxation of the formulation.
An empirical study demonstrates that the proposed approach can achieve better clustering accuracy than state-of-the-art
methods when the data is high-rank, the percentage of  missing data is high, or the subspaces are similar.
}%


\HISTORY{}

\maketitle

%


\section{Introduction}
\label{sec:intro}
We consider a real-valued matrix $X \in \mathbb{R}^{\ambientdimension \times \numvectors}$ in which each column $X_1, X_2,
\ldots, X_\numvectors$ is assumed to lie near one of $\numsub$ unknown subspaces, $\mathcal{S}_i$ with dimension
$\lowrank[i]< \ambientdimension$, $i=1,\ldots,\numsub$.
Given a set of $\observedset$ of observed entries of the matrix $X$,
\emph{subspace clustering with missing data} (SCMD) is the task of identifying clusters of the columns of
$X_{\observedset}$ belonging to the same subspace and inferring the subspace associated with each cluster of columns.

If the clustering of columns is known, then the subspace associated with each cluster can be estimated by applying
methods for low-rank matrix completion (LRMC)  on the columns in each cluster \citep{recht_2009, lrmc-near-optimal, lrmc-svd, laura-grouse,
lrmc-recht-bounds,lrmc-sampling-daniel,lrmc-review}. 
The SCMD problem has applications in machine learning for image classification \citep{MNIST, application-image-classification}, motion segmentation \citep{Hopkins155, motion-segmentation-vidal-2010}, and recommendation systems \citep{matrix-completion-recommendation}.

In this work, we propose a novel mixed-integer linear programming (MILP) solution framework for the SCMD problem that is
based on dynamically determining a set of candidate subspaces and optimally assigning columns (also called data points) to the closest selected subspace.
We refer to our method as \ourmethod: Mixed Integer Subspace Selector with Dynamic Subspace Generation.
A key contribution of our approach is a method for identifying, in a rigorous manner, a suitable set of candidate subspaces to include in the formulation.
We cast this subspace generation problem as a nonlinear, nonconvex optimization problem and propose a gradient-based approximate solution approach.
We also use Benders decomposition to solve the linear programming (LP) relaxation of the MILP, which allows our
framework to accommodate both a large number of candidate subspaces and a large number of data points.
The model has the advantage of integrating the subspace generation and clustering in a single, unified optimization
framework without requiring any hyperparameter tuning when the number of subspaces and subspaces dimensions are known.
When the number and/or dimensions of the subspaces  are unknown, we propose add a regularization term to the objective
based on the effective dimension of the selected subspaces.
Our computational study demonstrates that the proposed method can achieve higher clustering accuracy than
state-of-the-art methods when the underlying matrix $X$ is high-rank, the percentage of missing data is high, or the
subspaces are similar to each other.

{\bf Prior Work.}
\cite{review-maryam-2021} provide a a review of subspace clustering methods in the case that all entries of $X$ are
observed.
Most of the methods for subspace clustering with missing data have been adapted from the methods initially proposed for fully-observed data.
The tightest known conditions for union of subspaces identifiability with missing data have been established by
\citet{PimentelAlarcn2016TheIR}, where the authors show that for ambient dimension $\ambientdimension$ and low
dimensional subspaces of dimension $\lowrank$, observing $O(\lowrank \ambientdimension)$ columns per subspace is both necessary and sufficient for subspace
identification.

The dominant approach in subspace clustering is based on a self-expressiveness property, originally proposed for fully-observed data by \citet{ssc-ehsan-2009}.  
Self-expressive methods learn a sparse representation of the data by solving an optimization problem of the form: 
\begin{eqnarray}
	C^*=\argmin \|X-XC\|^2_F+\lambda \rho(C) \qquad \text{s.t. diag}(C)=0.
	\label{eq:ssc}
\end{eqnarray}
Self-expressive methods have been studied extensively for different choices of $\rho(\cdot)$, e.g., $\ell_1,\ell_2$, and nuclear norm \citep{ssc-ehsan-2009,LRR,LSR,ssc-lrr-zhuang,ssc-trace-lasso,ssc-Vidal-2013,ssc-ensc-music,ssc-ensc,ssc-lrr}.
The value $C^*_{ij}$ can be interpreted as a link or connection between data points $i$ and $j$.
The segmentation is obtained by applying spectral clustering on a graph $G$ with adjacency matrix $A = |C^*| + |C^*|^T$,
which uses $k$-means cluster of the eigenvectors of the Laplacian of $G$ \citep{Ng01onspectral}.
Self-expressive methods for subspace clustering have been extended to the case of missing data. 
Let $\observedset_\datavector$ denote the set of observed components of vector $\datavector\in [\numvectors]$ and $I_{\observedset} \in \{0,1\}^{\ambientdimension \times \numvectors}$ be the indicator matrix of observed entries such that $[I_{\observedset}]_{i\datavector}=1$ if $(i,\datavector)\in \observedset$, and $0$ otherwise. 
Let $\circ$ denote the Hadamard product.
\cite{ssc_missing_icml2015} proposed to zero-fill the missing entries in $X$ to get $X_{ZF}$ and to solve~\eqref{eq:ssc} while restricting the loss to observed entries, i.e.,
\begin{alignat}{1}
  & X_{ZF} = \begin{cases}
    X_{ij},& \text{if }(i,\datavector)\in \observedset, \text{i.e., } X_{i\datavector} \text{ is observed}\\
    0, & \text{if } (i,\datavector)\notin \observedset, \text{i.e., } X_{i\datavector} \text{ is not observed}\\
  \end{cases} \\
  &C^*=\argmin \|I_{\Omega} \circ(X_{ZF}-X_{ZF}C)\|^2_F+\lambda \rho(C) \qquad \text{s.t. diag}(C)=0.\label{eq:ssc-missing}
\end{alignat}      

      \cite{ssc-md-theory-Tsakiris} and \cite{ssc-md-theroy-charles} studied the theoretical conditions under which the solution to~\eqref{eq:ssc-missing} is subspace-preserving, i.e., each data point is only connected to points lying in the same subspace.

Representative of a class of methods that alternate between subspace estimation and assignment,
\cite{ssc_missing_icml2015} proposed to apply a matrix completion algorithm to recover the missing entries in $X$ and then solve~\eqref{eq:ssc}.
This approach has been observed to fail when the data matrix is high-rank, i.e., $\sum_{i=1}^\numsub \lowrank[i] \approx \ambientdimension$.
\cite{Lane_2019_ICCV} proposed to alternate between subspace clustering and group wise low-rank matrix completion (gLRMC).
\cite{laura-k-grouse} use GROUSE \citep{laura-grouse} for subspace estimation and assign each point to the orthogonally closest subspace.
They use  probabilistic farther insertion for initializing $\numsub$ subspaces.
In all alternating methods, the subspace estimation process is often faulty when an estimated cluster has points from
different subspaces.
In an extensive empirical evaluation of existing SCMD algorithms, \cite{Lane_2019_ICCV} concluded that the zero-filled elastic net subspace clustering method \citep{ssc-ensc} when alternated with low-rank matrix completion showed the overall best performance. This method is referred to as Alt-PZF-EnSC+gLRMC.

Some methods  pose the two problems of subspace estimation and assignment in a joint optimization framework, often resulting in complex, nonconvex problems \citep{structured-ssc-vidal, hrmc-elhamifar-NIPS2016,mc-fan-2017}.
Matrix factorization approaches have also been adopted to SCMD \citep{daniel-scmd,daniels-gscc}.
Empirical experiments in \citep{Lane_2019_ICCV} showed that these methods are outperformed by the alternating methods in terms of clustering error.


The work most closely related to our approach are the MILP-based methods for subspace clustering with fully-observed data. 
\cite{floss} were the first to propose an MILP-based method for subspace clustering called Facility Location for Subspace Segmentation (FLoSS).
FLoSS generates the candidate subspaces at random and formulates the subspace clustering problem as an
MILP.  The goal of the MILP is to minimize the orthogonal distances of data points to candidate subspaces such that $\numsub$ subspaces are selected, and each vector is assigned to a selected subspace.
\cite{minimal-basis-facility} extended the FLoSS model to Minimal Basis FLoSS (MB-FLoSS), where the 
subspace hypothesis generation strategy is based on finding the minimal basis subspace representation for the data matrix and relies on Low Rank Representation (LRR)\citep{LRR}.
\citet{hu-milp-2015} proposed the concept of constrained subspace model.
They integrated the facility-based model with manifold and spatial regularity constraints to develop a constrained subspace modeling framework. 
The number of candidate subspaces is small ($\leq 50$) in their experiments.
The method becomes inefficient when the number of candidate subspaces is higher, and the approach heavily relies on the efficiency of initial candidate subspaces generated for which they use over segmentation in LRR \citep{LRR}.
In particular, they generate more subspaces than the ground truth
(e.g, $2\times\numsub$) with LRR, and then use MILP to select $\numsub$ of them.
None of these existing MILP-based approaches account for missing data or scale to instances with a large number of candidate subspaces.
Moreover, all of the approaches require that candidate subspaces are explicitly enumerated as an input to the model, and either rely on random sampling or other subspace clustering algorithms for generating candidate subspaces.
Hence, these methods are incapable of correcting themselves based on the clustering quality.  Our approach 
\ourmethod can handle a large set of candidate subspaces through the use of Benders decomposition and identifies new candidate subspaces dynamically through the use of column generation.
Casting SCMD as an MILP offers several other advantages.
The formulation can easily be extended to incorporate prior information about the data, such as vectors lying in the same or different subspaces and bounds on number of subspaces.

The paper begins with a description of the MILP formulation in Section~\ref{sec:ip-formulation}.
Section~\ref{sec:decomposition} discusses our decomposition approach to solve the model, and
Section~\ref{sec:experiments} presents experimental results that show the effectiveness of our framework.

 \section{MILP Formulations}
 \label{sec:ip-formulation}
We consider a real-valued matrix $X \in \mathbb{R}^{\ambientdimension \times \numvectors}$ whose columns are concentrated near a union of $\numsub$ subspaces with dimensions $\lowrank[1],
\lowrank[2],\dots,\lowrank[K]$.
For an integer $\numcols$, we denote $[\numcols] := \{1,2,\ldots,\numcols\}$. We denote data vector (column)
$\datavector \in [\numvectors]$ as $X_\datavector$.
We assume that we observe a subset $\Omega \subseteq [\ambientdimension] \times [\numvectors]$ of the entries of $X$, and
given this data the goal of subspace clustering with missing data (SCMD) is to identify the $K$ subspaces together with assignment of
data points to subspaces. This consequently leads to a clustering of points and a method for estimating missing
entries of $X$. 
In Section \ref{subsec:ip-known-dim}, we assume that subspaces dimension $\lowrank[1],\dots,\lowrank[\numsub]$ are known.
We relax this assumption in Section \ref{subsec:ip-unknown-dim} and let the model self-determine the subspaces dimensions with the help of a regularized objective.
The matrix $X$ is referred to as \emph{low-rank} when $\sum_{k=1}^\numsub \lowrank[k]\ll \min\{\ambientdimension,\numvectors\}$
and as \emph{high-rank} when $\sum_{k=1}^\numsub \lowrank[k]\approx \min\{\ambientdimension,\numvectors\}$.

Our approach is based on iteratively building a (potentially very-large) collection of candidate subspaces. MILP is
employed to simultaneously select the best set of $\numsub$ candidate subspaces and assign each column of $X$ to its
closest selected subspace.  For each candidate subspace $\column \in [\numcols]$, we let $\basis[\column] \in
\mathbb{R}^{\ambientdimension \times \lowrank[t]}$  be a basis for its column subspace, where $ \lowrank[t]$ is the
dimension of candidate subspace $\column$.
We define the distance of vector $\datavector \in [\numvectors]$ to a candidate subspace $\column \in [\numcols]$ as the
sum of the squared residuals on the observed entries:
 \begin{equation}
 	\label{eq:distcomp}
 	\distance  := \min_{v\in \mathbb{R}^r} \Big\{ \sum_{i: (i,\datavector)  \in \observedset} (X_{i\datavector} - (\basis[t] v)_i)^2 \Big\}. 
 \end{equation}
These quantities have a closed-form solution in terms of a simple projection operator \citep{laura-k-grouse}.  
In particular, let $\basis[\observedset{},j]$ denote the restriction of the subspace $\basis$ to the rows observed in
column $\datavector$, and define the projection operator $\projectionOperator
:=\basis[\observedset,j](\basis[\observedset,j]^T \basis[\observedset,j])^{-1}\basis[\observedset,j]^T$. Then the
squared residual $\distance$ can be obtained  as 
\begin{equation}
	\distance = \|X_{\observedset,\datavector}-\projectionOperator[(\basis_t)_{\observedset,\datavector}](X_{\observedset,\datavector})\|_2^2 . 
	\label{eq:residual-closed-form}
\end{equation}

For fully-observed data, this is a natural choice for cost function since its value is zero if vector $\datavector$ is 
in subspace $\column$. However, with missing data, the choice of cost function becomes less clear since zero residual on observed entries for $\distance$ does not necessarily imply that vector $j$ lies perfectly on subspace $\column$. 
\cite{laura-k-grouse} showed that for a given fully observed vector $X_\datavector \in \mathbb{R}^{\ambientdimension}$, if
\begin{equation}
  \|X_\datavector-\projectionOperator[\basis_0](X_\datavector)\|<\|X_\datavector-\projectionOperator[{\basis[\column]}](X_\datavector)\|\ \forall \column \in [\numcols]\symbol{92} \{0\},
  \label{eq:trueres}
\end{equation}
then, with high probability, 
for the same data vector $X_\datavector$ but now partially observed
on $\observedset$ (sampled uniformly from $[\ambientdimension]$ with replacement), if ``enough'' elements of the 
data point $X_j$ are observed (see \cite{laura-k-grouse} for details on quantification of ``enough''), then
\begin{equation}
	\|X_{\Omega,\datavector}-\projectionOperator[(\basis_0)_{\observedset,\datavector}](X_{\observedset, \datavector})\|<
	\|X_{\Omega,\datavector}-\projectionOperator[(\basis_t)_{\observedset,\datavector}](X_{\observedset, \datavector})\| \ \forall t\in [T]\symbol{92} \{0\}.
	\label{eq:observedres}
\end{equation}
This implies that, with high probability, subspace assignment based on \eqref{eq:trueres} is the same as the one based on \eqref{eq:observedres}. 
We refer reader to \cite{laura-k-grouse} for more details.
We also note that a different cost model could also be incorporated into our framework. 

Next, we describe a model based on selecting $\numsub$ subspaces from a given collection of $[\numcols]$ subspaces for both known and unknown subspaces dimensions in Section \ref{subsec:ip-known-dim} and \ref{subsec:ip-unknown-dim} respectively.
In Section \ref{subsec:rowGeneration}, we discuss how to solve the  proposed model for a fixed set of subspaces using Benders decomposition.
This allows to solve the model efficiently for large $\numvectors$ and $\numcols$.
In Section \ref{subsec:columnGeneration}, we discuss how to generate new candidate subspaces dynamically with a column generation approach.
We finally discuss our unified framework \ourmethod in Section \ref{subsec:miss-dsg}.

Let $\assignvar \in\{0,1\}, \forall \datavector \in [\numvectors], \column \in [\numcols]$ be a binary assignment variable that takes value 1 if vector $\datavector$ is assigned to subspace $\column
 $, and $\selectvar \in \{0,1\}, \forall \column \in[\numcols]$ be a binary selection variable with $z_t=1$ if subspace $\column$ is selected.
 The assignment of points to selected subspaces is similar to the facility location problem, where the goal is to select
 which facilities to open and to assign each customer to an open facility. In our SCMD formulation, subspaces
 play the role of facilities, and vectors play the role of customers. Previously proposed MILP methods in
 the literature have this same facility-location structure \citep{floss,minimal-basis-facility,hu-milp-2015}. However, our model has some key differences:
 \begin{enumerate}[label=(\alph*)]
 	\item We allow missing data while existing MILP approaches restrict to fully observed data.
 	\item Our framework generates subspaces dynamically while existing approaches are heavily dependent on initially generating candidate subsspaces.
 	\item Our framework is capable of handling a larger number of candidate subspaces than existing methods through the use of Benders decomposition.
 \end{enumerate}
For improved readability, in Section \ref{subsec:ip-known-dim} we first discuss the formulation for the simpler case
where subspaces dimensions are assumed to be known and equal, i.e.,
$\lowrank[1]=\lowrank[2]=\cdots=\lowrank[k]=\lowrank$. We then relax this assumption in Section
\ref{subsec:ip-unknown-dim}. 

  \subsection{Known subspaces dimension }
 \label{subsec:ip-known-dim}
  Given $\numcols$ candidate subspaces each of dimension $\lowrank$, we formulate the SCMD problem as the MILP 
 \citep{floss,hu-milp-2015}:
\begin{subequations}
	\label{eq:milp} 
	\begin{alignat}{2}
		\min_{\assignvar[\null] \in \{0,1\}^{\numvectors \times \numcols},\selectvar[\null] \in \{0,1\}^\numcols} \ &  \sum_{\column\in [\numcols]} \sum_{\datavector \in [\numvectors]} \distance \assignvar \owntag[milp-obj]{MILP} \\
		& \sum_{\column \in [\numcols]} \assignvar = 1, && \forall \datavector \in [\numvectors]  \label{eq:assignment}\\
		& \assignvar \leq \selectvar, && \forall \datavector \in [\numvectors], \column \in[\numcols] \label{eq:facopen}\\
		&	\sum_{\column \in[\numcols]}\selectvar = \numsub. &&  \label{eq:maxnumz} 
	\end{alignat}
\end{subequations}
The objective \eqref{eq:milp-obj} seeks the least cost assignment of vectors to subspaces. 
Constraints \eqref{eq:assignment} ensures that each vector is assigned to exactly one subspace, and constraints
\eqref{eq:facopen} enforce that a vector can only be assigned to a selected subspace. 
Constraint \eqref{eq:maxnumz} requires that exactly $\numsub$ subspaces are selected.
The major difference in~\eqref{eq:milp} compared to  the models proposed in \citep{floss, hu-milp-2015} is that the
distance metric $\distance$ in~\eqref{eq:milp} is based on partial assignment cost~\eqref{eq:residual-closed-form}, as
we do not assume data is fully observed.

 \subsection{Unknown Subspaces Dimension}
 \label{subsec:ip-unknown-dim}
 A common occurrence in SCMD problems is that the number of subspaces and their dimension are unknown.
 The objective in model~\eqref{eq:milp} may be in appropriate in this case, because 
 it is likely to favor subspaces of higher dimension, since these are inherently more flexible and hence will lead to lower residuals.
 Thus, when the subspace dimensions are unknown, we propose to augment the objective in \eqref{eq:milp} with a complexity measure of the
 candidate subspaces, the effective dimension (\ed).
 \cite{gpca-mix-UoS}  defined effective dimension for $X$ on a union of subspace models $\textit{S}=\cup_{k=1}^\numsub S_k$ as follows:
 \begin{equation}
   \ed(X,S) := \frac{1}{\numvectors} \sum_{k=1}^\numsub \lowrank[k] (\ambientdimension-\lowrank[k])+\frac1n
	\sum_{k=1}^\numsub \numvectors_k\lowrank[k] .
   \label{eq:ed}	   
 \end{equation}
 The first term in the definition of \ed\ in \eqref{eq:ed}, $ \lowrank[k] (\ambientdimension-\lowrank[k])$ is the complexity of $\basis[\column]$---the number of real numbers needed to specify a $k$ dimensional subspace $S_k$ in $\mathbb{R}^\ambientdimension$.
 The second term of \eqref{eq:ed}, $\numvectors_k\lowrank[k]$ is the number of real numbers needed to specify the $\lowrank[k]$ coordinates of the $\numvectors_k$ sample points in the subspace $S_k$.

To allow for the model to trade-off accuracy with complexity, we add the \ed\ term with a weight parameter $\lambda$ in our objective as follows:
 \begin{subequations}
   \label{eq:milp-ed} 
   \begin{alignat}{2}
     \min_{\assignvar[\null] \in \{0,1\}^{\numvectors \times \numcols}, \selectvar[\null] \in \{0,1\}^{\numcols}} & \sum_{\column\in [\numcols]}\sum_{\datavector \in [\numvectors]} \distance \assignvar+ \frac{\regularizer}{\numvectors}\sum_{\column \in \numcols}\Big(\lowrank[t](\ambientdimension-\lowrank[t]) \selectvar+\sum_{\datavector \in [\numvectors]}\lowrank[t]\assignvar\Big) \\
     & \sum_{\column \in [\numcols]} \assignvar = 1, && \forall \datavector \in [\numvectors]  \label{eq:assignment-ed}\\
     & \assignvar \leq \selectvar, && \forall \datavector \in [\numvectors], \column \in[\numcols] \label{eq:facopen-ed}\\
     &	\sum_{\column \in[\numcols]}\selectvar = \numsub.  \label{eq:maxnumz-ed} 
   \end{alignat}
 \end{subequations}
Formulation~\eqref{eq:milp-ed} can handle subspaces of multiple dimensions and choose the best union of subspaces model by self-determining dimensions of subspaces.
Constraint \eqref{eq:maxnumz-ed} is included but can be removed if $\numsub$ is unknown. 
An important consideration in model~\eqref{eq:milp-ed} is the choice of regularization parameter $\regularizer$ which
accounts for the trade-off between lower assignment cost and complexity of the selected subspaces. A smaller value of
$\regularizer$ would promote model~\eqref{eq:milp-ed} to select higher complexity subspaces (basis with higher
dimensions $\lowrank[t]$) while a larger value of $\regularizer$ would promote the model to select lower complexity subspaces (basis with lower dimensions $\lowrank[t]$). We discuss this in more detail in Section \ref{subsec:penalty-choice}.
Note that~\eqref{eq:milp} is a special case of~\eqref{eq:milp-ed} with $\regularizer=0$.

\section{Decomposition Algorithm}
\label{sec:decomposition}

We next discuss how to solve formulation~\eqref{eq:milp-ed} for a given set of candidate subpsaces via Benders
decomposition and how to dynamically generate a set of candidate subspaces
to use in the formulation. 

 
As is standard for solving a MILP, the first step in solving formulation~\eqref{eq:milp-ed} is to solve its LP relaxation.
The LP relaxation of~\eqref{eq:milp-ed} is the problem created by replacing the integrality conditions $\selectvar, \assignvar \in \{0,1\}$ with simple bound constraints $\selectvar, \assignvar \in [0,1]$.
The optimal solution value of the LP relaxation provides a lower bound on the optimal solution to \eqref{eq:milp-ed}.
The number of candidate subspaces $\numcols$ and the number of points $\numvectors$ may be quite large, so solving the LP relaxation could be a computational challenge.
We discuss in Section~\ref{subsec:rowGeneration} how Benders decomposition can address this challenge for solving the LP
relaxation of \eqref{eq:milp-ed}.
In Section~\ref{subsec:columnGeneration}, we describe our approach for dynamically generating additional candidate
subspaces to include in the formulation, which is also done for the LP relaxation.
Finally, in Section~\ref{subsec:miss-dsg} we describe how these components are integrated into an overall approach for yielding
feasible solutions to \eqref{eq:milp-ed}.

\subsection{Benders Decomposition}
\label{subsec:rowGeneration}
Benders decomposition is a well-known technique to solve large LP problems that have special structure \citep{bnnobrs1962partitioning}.
It has been applied to large-scale facility locations by \citet{ufl-benders}. Since our SCMD
formulation~\eqref{eq:milp-ed} has the same structure, we can apply the same approach for solving its LP relaxation. 
The first step in the decomposition approach is to write a reformulation that eliminates the $\assignvar$ variables and
adds a continuous variable $\costvar$ for each $\datavector \in[\numvectors]$ that represents the assignment cost for
vector $\datavector$. The resulting reformulation of the LP relaxation of \eqref{eq:milp-ed} is 

 \begin{subequations}
	\label{eq:master-ed} 
	\begin{alignat}{2}
		\min_{\costvar[\null] \in R^n , \selectvar[\null]\in [0,1]^\numcols} & \sum_{j \in [n]} \costvar +
		\frac{\regularizer}{\numvectors}\sum_{\column \in \numcols}\lowrank[t](\ambientdimension-\lowrank[t])\selectvar\\
		& \costvar \geq \Phi_\datavector(\selectvar[\null]), && \forall \datavector \in [\numvectors]\\
		&	\sum_{\column \in[\numcols]}\selectvar = \numsub,  
	\end{alignat}
\end{subequations}
where $\Phi_j(\cdot)$ is defined as
\begin{equation}
  \Phi_\datavector(\hat{\selectvar[\null]})=\min_{\assignvar[\null]} \Big\{ \sum_{\column \in[\numcols]}(\distance+\frac{\regularizer}{\numvectors}\lowrank[t]) \assignvar[\column]: \sum_{\column \in[\numcols]} \assignvar[t] = 1, 0 \leq \assignvar[t] \leq \hat{\selectvar[\null]}_\column,  \forall \column \in[\numcols] \Big\}.
	\label{eq:phi}
\end{equation}
The function $\Phi_\datavector(\hat{\selectvar[\null]})$ computes the minimum (fractional) assignment cost for the vector $\datavector \in[\numvectors]$ for a given
(possibly fractional) vector of facility opening decisions $\hat{\selectvar[\null]} \in [0,1]^\numcols$.
The function $\Phi_\datavector(\cdot)$ is piecewise-linear and convex,
and Benders decomposition works by dynamically building a lower-bound approximation to $\Phi_\datavector(\cdot)$ via the
addition of Benders cuts.

To simplify notation in the derivation of the Benders cuts, we define $\cost =
\distance+\frac{\regularizer}{\numvectors}\lowrank[t]$ for $j \in [n]$ and $t \in [\numcols]$.
The optimization problem \eqref{eq:phi} used to evaluate $\Phi_\datavector(\hat{\selectvar[\null]})$ has a closed-form solution.
Moreover, its evaluation also gives sufficient information to derive the Benders cuts that define the lower-bounding approximation.
For each $j \in [n]$, let $\{\sigma^\datavector_1,\ldots,\sigma^\datavector_\numcols\}$ be a permutation of $[\numcols]$ satisfying $c_{\datavector\sigma^\datavector_1} \leq c_{\datavector\sigma^\datavector_2} \leq \cdots \leq \cost [\datavector \sigma^\datavector_\numcols]$, and let $\column^*_\datavector := \min\{ \column: \sum_{s=1}^\column \hat{\selectvar[\null]}_{\sigma^\datavector_s} \geq 1 \}$ be the \emph{critical index}.
In other words, the \textit{critical index} is the index of the costliest subspace to which any portion of vector $\datavector$ is assigned.
As described in \cite{ufl-benders}, the Benders cut that can be used to lower-approximate the function $\Phi_j(\cdot)$ is
\begin{equation}
  \costvar+\sum_{i=1}^{\column_\datavector^*-1}(\cost[\datavector\sigma^\datavector_{\column_\datavector^*}]-\cost[\datavector\sigma^\datavector_i])\selectvar[\sigma^\datavector_i]\geq \cost[j\sigma^j_{t_j^*}]. \label{eq:benders} 
\end{equation}

These inequalities are accumulated iteratively.  
Let $p_\datavector$ denote the number of Benders cuts included in the model at the current stage in the algorithm for each $\datavector \in[\numvectors]$.
Let $\column^*_{\datavector\cind}$ denote the critical index for vector $\datavector \in[\numvectors]$ associated with Benders cut $\cind\in[p_\datavector]$, and let $\cost[\datavector \cind]^* = c_{j\sigma^j_{t^*_{j\cind}} }$ denote the critical cost for the $\datavector^{th}$ vector in cut $\cind\in[p_\datavector]$. 
The Benders \emph{master problem} is then
 \begin{subequations}
\label{eq:master-lp-relaxation}
\begin{alignat}{2}
  \min_{\costvar[\null] , \selectvar[\null]} & \sum_{\datavector \in [\numvectors]} \costvar +\frac{\regularizer}{\numvectors}\sum_{\column \in [\numcols]}\lowrank[\column](\ambientdimension-\lowrank[\column])\selectvar   \\
  & \costvar+\sum_{\ell=1}^{\column^*_{j\cind}-1}(\cost[\datavector\cind]^* -c_{j\sigma^i_\ell})z_{\sigma^j_\ell}\geq c^*_{j\cind} , && \ \forall \datavector \in[\numvectors], \cind \in [p_\datavector],  \owntag[alpha-constr]{$\alpha_{\datavector\cind}$} \\
  &\sum_{\column \in[\numcols]}\selectvar = \numsub, \owntag[beta-constr]{$\beta$} \\
  & 0\leq \selectvar \leq 1, && \forall  \column \in[\numcols].\owntag[mu-constr]{$\mu_t$} 
\end{alignat}
 \end{subequations}
Here $\alpha,\beta$ and $ \mu$ are dual variables corresponding to the respective constraints, and will play an important role in the column generation process described in Section~\ref{subsec:columnGeneration}.
For $\numcols>\numsub$, LP~\eqref{eq:master-lp-relaxation} is feasible and bounded, and hence an optimal solution $(\hat{\costvar[{\null}]},\hat{\selectvar[{\null}]})$ exists.
The subproblem \eqref{eq:phi} is solved to evaluate $\Phi_\datavector(\hat{\selectvar[{\null}]})$ for each $\datavector \in [\numvectors]$, and to generate new Benders cuts \eqref{eq:benders}.  
If $\Phi_\datavector(\hat{\selectvar[{\null}]}) = \hat{\costvar[\null]}_j$, then the generated inequality does not improve the approximation to $\Phi_\datavector(\cdot)$, and the cut is not added to \eqref{eq:master-lp-relaxation}.
The Benders procedure stops when no new cuts are added.
At this point, the LP relaxation of \eqref{eq:milp-ed} is solved.

\subsection{Column Generation}
\label{subsec:columnGeneration}
In our discussion to this point, we have assumed that we are given $\numcols$ candidate subspaces.
Key to our approach is a \emph{column generation} method for dynamically identifying new subspaces that have the potential to improve the solution to \eqref{eq:milp-ed}.
Column generation is a classical method for solving large-scale LP \citep{column-generation-ford} that also has seen
significant use in solving MILP problems \citep{barnhart.et.al:98}. We apply column generation to the LP relaxation of
\eqref{eq:milp-ed}, or more specifically, to the Benders reformulation of this LP relaxation 
\eqref{eq:master-lp-relaxation}.

The key idea behind column generation is to create an auxiliary problem, called the \emph{pricing problem}, whose solution identifies if there is an additional variable (candidate subspace), that, when added to the LP \eqref{eq:master-lp-relaxation}, could improve its solution value.
The formulation of the pricing problem follows from LP duality theory.
If the \emph{reduced cost} of a column (subspace variable) is negative, then, by increasing the value of that variable from its nominal value of zero, the objective value of the LP may decrease.
Thus, we seek columns (subspaces) with negative reduced cost.
If all columns have non-negative reduced cost, the current solution of the LP with the set $[\numcols]$ of candidate
subspaces is optimal and adding columns to $[T]$ would not decrease the LP relaxation value.

Consider an arbitrary subspace variable $\selectvar$. Given the optimal dual variables ($\hat{\alpha},\hat{\beta}$) to the solution  of \eqref{eq:master-lp-relaxation}, the reduced cost of a variable $\selectvar$ is given by the formula
\begin{equation}
	\frac{\regularizer}{\numvectors}\lowrank[\column](\ambientdimension-\lowrank[\column]) - \sum_{\datavector \in [\numvectors]} \sum_{\cind \in [p_\datavector]} \hat{\alpha}_{\datavector\cind}\max\{{\cost[\datavector\cind]^*-\cost}, 0\} -\hat{\beta},
	\label{eq:reduced-cost}
\end{equation}
where the term $\max\{{\cost[\datavector\cind]^*-\cost}, 0\}$ that is used to record the contribution of Benders cut $\cind
\in [p_\datavector]$ for data vector $\datavector \in [n]$ is used to reflect the fact that the coefficient
$\cost[\datavector\cind]^*-\cost$ only appears for variable $\selectvar$ in this cut when $\cost[\datavector\cind]^* > \cost$.
	
We formulate the problem of finding a subspace that corresponds to a column having negative reduced cost as a problem of
finding a basis matrix
$U \in \mathbb{R}^{\ambientdimension \times \lowrank}$, where $\lowrank$ is the candidate dimension of the subspace.   
In order to derive the reduced cost for the variable associated with a subspace defined by  basis matrix $U$, we recall
from \eqref{eq:distcomp} that for a subspace defined by basis matrix $U_t$,
\begin{equation}
	\label{eq:hi}
	\cost = \distance +  \frac{\regularizer}{\numvectors} \lowrank = \min_{v\in \mathbb{R}^{\lowrank}} \Big\{ \sum_{i: (i,\datavector)  \in \observedset} (X_{i\datavector} - (\basis[t] v)_i)^2 \Big\} 
	+ \frac{\regularizer}{\numvectors} \lowrank =: h_\datavector(\basis[\column]) .
\end{equation}
Thus, substituting the function $h_\datavector(\basis)$ in for $\cost$ in the reduced cost expression
\eqref{eq:reduced-cost}, the problem of finding a column corresponding to a rank-$\lowrank$ basis of minimum reduced
cost can be formulated as:
\begin{equation}
	\min_{\basis \in \mathbb{R}^{\ambientdimension \times \lowrank}} g_r(\basis) :=
 \frac{\regularizer}{\numvectors} \lowrank(\ambientdimension-\lowrank)  - \sum_{\datavector \in [\numvectors]}
 \sum_{\cind\in [p_\datavector]} \hat{\alpha}_{j\cind} 
	\max\{ \cost[\datavector\cind]^* - h_\datavector(\basis) , 0 \} -\hat{\beta} .
	\label{eq:pricing}
\end{equation}
The objective in problem~\eqref{eq:pricing} is not convex, and hence we find locally minimal solutions
to~\eqref{eq:pricing} with a gradient-based method. To explore subspaces of varying dimensions, we
solve~\eqref{eq:pricing} for $\lowrank \in \{1,2,\dots,\lowrank[max]\} $ where $\lowrank[max]$ is an upper bound on the
subspace dimension, which we assume is given as an input to the model.

{\bf Gradient-based method for solving \eqref{eq:pricing}. }
Consider a fixed $\lowrank\in \{1,2,\dots,\lowrank[max]\}$.
Observe that if $h_j(\basis) \neq \cost[ji]^* \ \forall \datavector \in [\numvectors], i \in [p_\datavector]$, then the
function $g_{\lowrank}(\cdot)$ is differentiable at $\basis$.
The partial derivative of $g_{\lowrank}(\cdot)$ with respect to matrix element $\basis[ab]$ evaluated at current iterate $\hat{U}$ is given by
\begin{equation}
  \frac{\partial g_{\lowrank}(\hat{\basis})}{\partial \basis_{ab}} = -\sum_{\datavector \in [\numvectors]} \sum_{\substack{\cind\in[p_\datavector] :\\ \cost[\datavector\cind]^*-h_\datavector(\hat{\basis})>0}} 
  2 \hat{\alpha}_{\datavector \cind}\sum_{\ell\in \Omega_\datavector} (X_{\ell \datavector}-\hat{u}_\ell^\top \hat{v}_\datavector)\hat{v}_{\datavector b} \quad \forall a\in [\ambientdimension],b\in[\lowrank],
  \label{eq:g_derivative}
\end{equation}
where $\hat{u}_\ell$ represents row $\ell$ of $\hat{U}$ and $\hat{v}_\datavector$ is the optimal solution to problem \eqref{eq:hi} that is solved when evaluating
$h_\datavector(\hat{\basis})$ for each $\datavector \in [\numvectors]$.
By adding the condition  $\cost[\datavector\cind]^*-h_\datavector(\hat{\basis})>0$ in the summation, we implicitly use a subgradient contribution of 0 at points of non-differentiability. 
We denote by $\nabla g_{\lowrank} (\hat{U})$ the $\ambientdimension \times \lowrank$ matrix of partial derivatives of
$g_{\lowrank}(\cdot)$ at $\hat{U}$.
\begin{algorithm}[h]
  \caption{Locally solving pricing problem for fixed subspace dimension}
  \label{alg:pricing}
  \SetAlgoLined
  \KwData{ $X_{\observedset}$}
  \KwInput{$\basis[0]\in R^{\ambientdimension \times \lowrank}$, maxIt=500, $\epsilon=0.001$ 
    \label{alg:pricing-input}
    \tcc*{initial subspace}}
  
  $\iter =0, \hat{\basis}=\basis[0], \mathcal{\basis}=\{\}$ \tcc*{iteration count}
  \While{ $\iter < \operatorname{maxIt}$ {\tt and} $\|\nabla g_r(\hat{\basis})\|>\epsilon$}{
    \For{$\datavector = 1,2,\dots,\numvectors$ \label{alg:pricing-projection-start} }{
      $\hat{v}_j=(\hat{\basis}_{\observedset,\datavector}^\top \hat{\basis}_{\observedset,\datavector})^{-1}
      \hat{\basis}_{\observedset,\datavector}(X_j)_\Omega$ \label{alg:projcalc}
      
      $\cost[\datavector, \iter]= \|(X_\datavector)_{\observedset}\; -\hat{\basis}_{\observedset,\datavector}\hat{v}_\datavector\|_2^2 + \frac{\regularizer}{\numvectors} \lowrank$ 
    }
    \label{alg:pricing-projection-end}
    
    Calculate $\nabla g(\hat{\basis})$ using \eqref{eq:g_derivative} 		\label{alg:pricing-gradient}
    \tcc*{Requires $\hat{\basis},\hat{v}_\datavector \ \forall \datavector\in [\numvectors]$} 
    
    Calculate $\tilde{g}$ using~\eqref{eq:g-approx}
    \label{alg:pricing-g-approx}
    
    $\gamma \gets \min(0.1,\frac{\tilde{g} -g(\hat{\basis})}{\|\nabla g(\hat{\basis})\|_F^2})$
    \label{alg:pricing-polyak}  \tcc*{Polyak step size}
    
    $\hat{\basis} \gets \hat{\basis} -\gamma\nabla g(\hat{\basis})$ \label{alg:pricing-grad-step} \tcc*{gradient step}
    
    $\iter \gets \iter +1$
    \If{$g_r(\hat{\basis}) < 0$}{
    $ \mathcal{\basis}\gets  \mathcal{\basis} \cup \hat{\basis}$
    \label{alg:pricing-store-U}
    }
  }
  \KwOutput{Set of subspace bases $\mathcal{\basis}$}
\end{algorithm}

We outline our gradient-based approach for locally solving pricing problem~\eqref{eq:pricing} for a fixed subspace dimension $\lowrank$ in Algorithm \ref{alg:pricing}. 
Since this method is not guaranteed to find a global optimal solution to the nonconvex problem~\eqref{eq:pricing}, we
run the method multiple times with different random choices of $\basis[0]$ to identify different locally-optimal solutions.
Our method for randomly initializing $\basis[0]$ is motivated by the fact that when the matrix is fully observed,
$\lowrank+1$ vectors per subspace are necessary and sufficient for subspace clustering. 
Hence, we randomly sample N$(>\lowrank+1)$ vectors from the $M$ vectors in the current LP solution of
\eqref{eq:master-lp-relaxation} with the largest $\hat{w}_j$ values, where $M > N$.  
In our experiments we use $M=5\lowrank[max]$ and $N=2\lowrank$.
This approach selects a subset of vectors with high residuals in the current solution and initializes the gradient descent algorithm with a best-fit subspace on that subset of vectors.
Then, we use a fast low-rank matrix completion algorithm, GROUSE \citep{laura-grouse}, to find the basis
$\basis[0]$ for a best-fit subspace for the sampled vectors. 
This $\basis[0]$ is provided as an input to the Algorithm \ref{alg:pricing} (line~\ref{alg:pricing-input}).

To calculate the gradient using equation (\ref{eq:g_derivative}), we first need to 
solve the optimization model in \eqref{eq:hi}, with $U_t = \hat{U}$, for each $\datavector \in [\numvectors]$.
The solution to this problem can be computed using the matrix
$(\hat{\basis}_{\observedset,\datavector}^\top \hat{\basis}_{\observedset,\datavector})^{-1}\hat{\basis}_{\observedset,\datavector}$
as used in line \ref{alg:projcalc} \citep{laura-k-grouse}.  

An important choice for the practical performance of gradient-based methods is the \textit{step size}.
We use the Polyak step size \citep{polyak}.
The Polyak step size formula requires the optimal value of objective function, $g^*$.
Since the optimal value $g^*$ is unknown, we approximate it with $\tilde{g}$ (in line~\ref{alg:pricing-g-approx}) as follows 
\begin{equation}
\tilde{g}\approx \frac{\regularizer}{\numvectors}\lowrank(\ambientdimension-\lowrank)-\sum\limits_{\datavector\in [\numvectors]}
\sum\limits_{\substack{\cind\in [p_\datavector]:\\	\cost[\datavector\cind]^*-h_\datavector(\hat{\basis})>0}}
\hat{\alpha}_{\datavector\cind}
c^*_{\datavector\cind}.
\label{eq:g-approx}
\end{equation}
This choice assures that $\tilde{g}$ is a lower bound on $g^*$, since it replaces $h_j(\basis)$ in the objective of
\eqref{eq:pricing} with $0$. 
The Polyak step size is then calculated as ${\gamma} \gets \frac{\tilde{g} -g(\hat{\basis})}{\|\nabla g\|_2^2}$.
To prevent the step size from getting too large when $\|\nabla(g)\|$ becomes small, we set  ${\gamma} \gets \min(0.1,\frac{\tilde{g} -g(\hat{\basis})}{\|\nabla g\|_2^2})$ (line~\ref{alg:pricing-polyak}).
These choices were made based on preliminary empirical experiments. We discuss this in more detail in Section \ref{subsec:polyak}. 
The gradient step is described in line~\ref{alg:pricing-grad-step}.
We terminate when the norm of the gradient becomes sufficiently small  ($\|\nabla g (\hat{U})\|\leq \epsilon = 0.001$ in
our implementation) or we reach the maximum number of allowed iterations (500 in our implementation). 
As described in line~\ref{alg:pricing-store-U}, we store the basis $\hat{U}$ in each iteration if it has negative
reduced cost, since each 
defines a potential candidate subspace to be added to the master problem~\eqref{eq:master-lp-relaxation}.

\newcommand{\rnc}{\operatorname{root\_node\_continue}}
\newcommand{\gencuts}{\operatorname{generate\_cuts}}
\subsection{MISS-DSG: Mixed Integer Subspace Selector with Dynamic Subspace Generation}
\label{subsec:miss-dsg}
\begin{algorithm}[]
  \SetAlgoLined
  \KwData{$X_\Omega$}
  \KwInput{max dimension $\lowrank[max]$, min and max \# of multi-starts: $\eta_{\min}, \eta_{max}$, max \# iterations $i_{\max}$  }
  \KwOutput{Partition of $[\numvectors]$ into $\numsub$ clusters: $S_k$, $k=1,\ldots,K$}
  Initialize MILP model \eqref{eq:milp-ed} with $[\numcols]$ as $\numinitbasis$ random subspaces of dimension $\lowrank[i]$
  for each $i\in [\lowrank[max]]$  \label{alg:MILP-init} \;
  Calculate $\cost$ for $\datavector\in [\numvectors],\column\in [\numcols]$ \label{alg:MILP-init-cost} \;
  $\rnc \gets {\tt True}$, \text{generate cuts} $\gets {\tt True}$, it $\gets 0$\;
  \While{$\rnc$ and $\iter < i_{\max}$\label{alg:MILP-root-node-loop-start}}{
    $\rnc \gets {\tt False}$ \tcc*{switched back on if new columns found}
    $\iter \gets \iter+1$\;
    \tcp{generate Benders cuts}
	 $\gencuts \gets {\tt True}$\; 
    \While{$\gencuts$ \label{alg:MILP-cutloopstart}
    }
	 {
      $\gencuts \gets {\tt False}$ \label{alg:MILP-benders-loop-start}\;
      solve master LP relaxation \eqref{eq:master-lp-relaxation} to obtain ($\hat{\costvar[\null]},\hat{\selectvar[\null]}$)\; \label{alg:MILP-lprelax}
      \For{$j = 1,2,\dots,\numvectors$ \label{alg:MILP-bendersstart}}{
        \uIf{ $ \hat{\costvar[\null]}_j<\Phi_\datavector(\hat{\selectvar[{\null}]})$\label{alg:MILP-cut-condition}}{
          Add Benders cuts of the form \eqref{eq:benders} to master \eqref{eq:master-lp-relaxation}\;
         $\gencuts \gets {\tt True}$
        } \label{alg:MILP-checkcuts}
      }\label{alg:MILP-bendersend}
    }    \label{alg:MILP-cutloopend} 
    \tcp{generate new columns}
    \For{$\lowrank=1,\dots,\lowrank[max]$
      \label{alg:MILP-startCG}		
    }
    {
      \For{$\eta=1,\dots,\eta_{\max}$
        \label{alg:MILP-multi-start-loop}}{
        Initialize $\basis[0]$ with a best-fit subspace for $2r$ randomly-sampled  vectors from the $5\lowrank[max]$ vectors with largest $\hat{\costvar[\null]}_\datavector$
        \label{alg:MILP-U0}\; 
        $\mathcal{\basis}\gets$Solve pricing problem using Algorithm \ref{alg:pricing} to generate candidate subspaces \label{alg:MILP-solve-pricing}\;
        \label{alg:MILP-lookup-red-cost-cols}
        \uIf{$\mathcal{\basis}\neq \emptyset$  \label{alg:MILP-neg-red-start}}{
          Extend $[\numcols]$ to include a new $\selectvar$ variables for each subspace $U_t \in \mathcal{\basis}$\;
          \label{alg:MILP-addNewCol}
          \uIf{$\eta>\eta_{\min}$\label{alg:MILP-breakCG-start}}{
            $\rnc \gets True$\; 
            \textbf{break} \tcc*{New columns found and minimum multi-starts done}
          }\label{alg:MILP-breakCG-end}
        }
      }
    } \label{alg:MILP-endCG}
    \uIf{$\operatorname{root\_node\_continue}$}{
      remove all Benders cuts from \eqref{eq:master-lp-relaxation}
      \tcc*{invalid due to new $\selectvar$ vars}
      \label{alg:MILP-delBenders}
    }
  } \label{alg:MILP-root-node-loop-end} 
  $\hat{x}_{\datavector\column}, \hat{z}_{\column} \gets$ Solve MILP model \eqref{eq:master-ed} with
  $\selectvar\in\{0,1\}^\numcols$, using a callback routine to search for violated Benders cuts at integer feasible solutions.
  \label{alg:MILP-solveMIP}\; 
  \textbf{return} $\{S_\column = \{\datavector\in[\numvectors]: \hat{\assignvar[\null]}_{\datavector \column}=1\}, \forall \column \in [\numcols] \text{ s.t. } \hat{\selectvar[\null]}_\column=1\}$ \; \label{alg:MILP-return} 
  \caption{Mixed Integer Subspace Selector with Dynamic Subspace Generation}
  \label{alg:miss-dsg}
\end{algorithm}

We now describe how the Benders decomposition and column generation processes are integrated into an overall solution
method, MISS-DSG, for the SCMD problem. 
 Algorithm \ref{alg:miss-dsg} describes the details.
 We initialize the model \eqref{eq:master-lp-relaxation} with $\numinitbasis$ randomly generated subspaces for each possible subspace dimension (line~\ref{alg:MILP-init}).
 In particular, to generate  a subspace of dimension $\lowrank$, we sample a matrix from a uniform distribution over
 $[-1,1]^{\ambientdimension\times \lowrank}$, and perform singular value decomposition (SVD) to get a basis. 
We then solve the master LP relaxation \eqref{eq:master-lp-relaxation}  in line~\ref{alg:MILP-lprelax}, and  generate Benders cuts for each $\datavector\in[\numvectors]$ (lines \ref{alg:MILP-bendersstart}-\ref{alg:MILP-bendersend}). 
Observe that if $\Phi_\datavector(\hat{\selectvar[{\null}]}) = \hat{\costvar[\null]}_j$, then the generated inequality does not improve the approximation to $\Phi_\datavector(\cdot)$, and the cut is not added to \eqref{eq:master-lp-relaxation} (line~\ref{alg:MILP-cut-condition}).
We repeat lines \ref{alg:MILP-bendersstart}-\ref{alg:MILP-bendersend} until no violated cuts are found.

We next proceed to generate new columns (lines~\ref{alg:MILP-startCG}-\ref{alg:MILP-endCG}) by solving the pricing problem~\eqref{eq:pricing} in order to look for negative reduced cost columns. 
We solve the pricing problem for every candidate dimension $\lowrank\in \{1,2,\dots,\lowrank[max]\}$ (line~\ref{alg:MILP-startCG}).
For each dimension $\lowrank$, we do multi-start (line~\ref{alg:MILP-multi-start-loop}) and locally solve the pricing
problem~\eqref{eq:pricing} using Algorithm \ref{alg:pricing} and store all bases having negative reduced cost in $\mathcal{\basis}$ (line~\ref{alg:MILP-solve-pricing}).
If any such basis are found, we add variables representing the possibility to select them to the master LP (line~\ref{alg:MILP-addNewCol}).
We proceed to the next dimension if any negative reduced cost
columns are found and we have exceeded the minimum number of multi-starts, $\eta_{\min}$ (lines~\ref{alg:MILP-breakCG-start}-\ref{alg:MILP-breakCG-end}).
We perform a maximum of $\eta_{\max}$ multi-starts for each dimension $\lowrank$. In our experiments, we use $\eta_{\min}=5$ and $\eta_{\max}=15$.

If new columns are found in the column generation process, we
delete existing Benders cuts from the master problem since they become invalid due to new $\selectvar$ variables  (line \ref{alg:MILP-delBenders}), and return to the process of generating Benders cuts. 
We repeat cut generation (lines~\ref{alg:MILP-cutloopstart}-\ref{alg:MILP-cutloopend}) and column generation
(lines~\ref{alg:MILP-startCG}-\ref{alg:MILP-endCG}) as long as we find new columns having negative reduced cost or until
we reach the maximum iteration limit, $i_{\max}$ which we set to be 15 in our experiments.

After we exit the root node loop we pass the updated MILP model with the new columns and cuts included to a MILP solver
(line~\ref{alg:MILP-solveMIP}) to be solved by a branch-and-cut method. 
We do not generate new columns ($\selectvar$ variables) after this point -- in other words, we generate columns only at the
root node, and not within the branch-and-cut phase of the solution of the MILP problem. As a result, even if we 
solve the pricing problem  \eqref{eq:pricing} to global optimality, this method is not guaranteed to find a globally
optimal solution due the possibility of missing needed columns at nodes in the branch-and-bound tree.
We do, however, generate Benders cuts within the branch-and-cut method for solving the MILP problem
as this is necessary to assure that integer feasible solutions encountered in the solution process have their costs correctly recorded before allowing them
to be accepted as an incumbent (best-known) solution.
This is implemented by defining a \lstinline[basicstyle=\ttfamily]|lazy constraint callback| that is called by the
solver any time it encounters an integer-feasible solution $(\hat{\costvar[\null]},\hat{\selectvar[\null]})$.
Within the lazy constraint callback, we check if $\hat{\costvar[\null]}_\datavector <
\Phi_\datavector(\hat{\selectvar[{\null}]})$ for any
$\datavector \in [\numvectors]$, and add the corresponding 
Benders cut of the form \eqref{eq:benders} to the model if so. The solver then adds these cuts to its formulation (and
consequently excludes the solution that had incorrect cost value $\hat{\costvar[\null]}_\datavector$).
After the MILP solver completes its branch-and-cut process,
we use the optimal solution returned to determine the selected subspaces and mapping of each vector to a selected subspace (line \ref{alg:MILP-return}). In particular, we select subspaces $\column \in [\numcols]$ such that $\hat{\selectvar[\null]}_\column=1$. 

\section{Effect of Algorithm Components}
\label{sec:modeling-choices}
We next report results of computational studies designed to investigate the importance of various components of MISS-DSG.

\subsection{Experimental Setup}
\label{subsec:experimental-setup}

We use Gurobi 8.1 for solving the LP and MILP problems and Python as the programming language. We set a time limit of 5000s  for each \ourmethod run. 
The computational study is conducted on a cluster of 4 core machines with a RAM of 16GB with Xeon X5690 CPU running at 3.46GHz.
We report results on instances generated randomly in a fashion similar to \citep{Lane_2019_ICCV}.  Specifically, we construct
$\numsub$ random subspaces with bases $\basis[k] \in \mathbb{R}^{\ambientdimension \times \lowrank[k]}, \forall
k\in[\numsub]$ by sampling entries from a standard Gaussian distribution. We then generate $\numvectors$ different data
vectors. Each data vector $\datavector \in[\numvectors]$ is sampled from one of the $\numsub$ subspaces, i.e.,
$X_\datavector=\basis[k]v_\datavector$ for a uniform random $k\in[\numsub]$ and $v_\datavector \in
\mathbb{R}^{\lowrank[k]}$ sampled from a standard Gaussian distribution. 
After generating data matrix $X$, we uniformly at random  drop a percentage $\missingpercentage$ of the entries in $X$
leaving the remaining entries as the observed components $\Omega$. 

\subsection{Impact of Benders Decomposition}
\label{subsec:bnb-comp}

We first demonstrate that Benders decomposition gives a very significant speedup in the time required to solve the LP in \ourmethod.
We compare the solution times for solving the LP relaxation without Benders decomposition (i.e., directly solving LP
relaxation of~\eqref{eq:milp-ed}) and with Benders decomposition  to solve the Benders formulation~\eqref{eq:master-lp-relaxation}. 
For this experiment, we generate an identical set of subspaces $([\numcols])$ to use in the two formulations and solve
the resulting LP relaxations, and we use $\regularizer=0$.
All results reported in this section are averaged over five different random instances for each combination of
parameters. 
 
In Table \ref{tab:benders-comp- num-cols}, we present the results for varying  number of candidate subspaces
($|\numcols|$), while fixing $\ambientdimension=30,\numvectors=200,\numsub=6,\missingpercentage=0,\lowrank[k] = 3\
\forall k\in[\numsub]$,
and in Table~\ref{tab:benders-comp-num-vectors}, we present the results for varying $\numvectors$ with fixed
$\ambientdimension=30, \numsub=6, \missingpercentage=0, |\numcols|=500, \lowrank[k]=3\ \forall k \in [\numsub]$. From
these two tables we observe that directly solving the LP relaxation of \eqref{eq:milp}  becomes very time-consuming as
either the number of candidate subspaces or the number of data points increases, whereas the time grows much more
modestly when using Benders decomposition. Thus, we find that the use of Benders decomposition is essential for the computational viability of
\ourmethod.

%
\begin{table}
	\begin{minipage}{0.5\linewidth}
	\centering
	\caption{Effect of $\mathbf{|\numcols|}$ on LP relaxation time (s) for $\mathbf{\ambientdimension=30,\numvectors=200,\numsub=6,\missingpercentage=0,\lowrank[k]=3\ \forall k \in [\numsub]}$}
	\label{tab:benders-comp- num-cols}
	\begin{tabular}{lrr}
		\toprule
		$|\numcols|$ & Without Benders & With Benders \\
		\midrule
		100 & 0.5 & 0.1\\
		500 & 3.2 & 0.3\\
		1000 & 48.2 & 1.6\\
		2000 & 729.8 & 2.8\\
		4000 & 1380.3 & 3.1\\
		\bottomrule
\end{tabular}
	\end{minipage}
\begin{minipage}{0.5\linewidth}
		\centering
		\caption{Effect of  $\mathbf{\numvectors}$ on LP relaxation time (s) for $\mathbf{\ambientdimension=30,\numsub=6,\missingpercentage=0,
			|\numcols|=500, \lowrank[k]=3\ \forall k\in [\numsub]}$ }
	\label{tab:benders-comp-num-vectors}
	\begin{tabular}{lrr}
		\toprule
		$\numvectors$ & Without Benders& With Benders\\
		\midrule
		100 & 7.4 & 0.5\\
		200 & 3.2 & 0.3\\
		400 & 109.5 & 0.8\\
		600 & 617.6 & 2.1\\
		1200 & 1147.3 & 4.5\\
		\bottomrule
	\end{tabular}
\end{minipage}
\end{table}

\subsection{Step-size Rule for Solving \eqref{eq:pricing}}
\label{subsec:polyak}
As discussed in Section~\ref{subsec:columnGeneration}, we propose to use the Polyak step size rule \citep{polyak} 
when solving the pricing problem~\eqref{eq:pricing}.
A major advantage of the Polyak rule is that it is does not require tuning the initial step size, as is required when
using a constant or diminishing step size rule.
We illustrate this on a test instance with $\lambda=0$, $\ambientdimension=20, \numvectors=200,\numsub=5, f=35\%,\lowrank[k]=4\
\forall k=1,\dots,4$.  
Figure \ref{fig:step-size-comp} displays the results, where we compare the evolution of the objective obtained using Algorithm \ref{alg:pricing} with the Polyak step size and with
a decaying step size ($\alpha_0/\texttt{it}$, where \texttt{it} is the iteration number). We plot the pricing problem
objective value~\eqref{eq:reduced-cost} on the y-axis and iteration number on x-axis. For the decaying step size rule,
we consider initial step size $\alpha_0\in \{0.001,0.01,0.1,1\}$. We observe that the Polyak step size leads to the
fastest convergence. We found similar behavior on other test instances.

 \subsection{The Necessity of Multi-start}
 We next discuss the importance of doing multi-start when solving the pricing problem~\eqref{eq:pricing}. We consider
 the same instance as discussed in Section~\ref{subsec:polyak} and solve the pricing problem with different
 initialization points. As shown in Figure \ref{fig:multi-start}, we observe that three different choices lead to three
 different local solutions. Hence, multi-starting can help find multiple different local minima and therefore improve
 the chances of identifying subspaces with negative reduced cost.
\begin{minipage}{\linewidth}
	\centering
	\begin{minipage}{0.45\linewidth}
		\begin{figure}[H]
			\caption{Comparison of Polyak step size with decay step size for same starting point}
			\includegraphics[width=\linewidth]{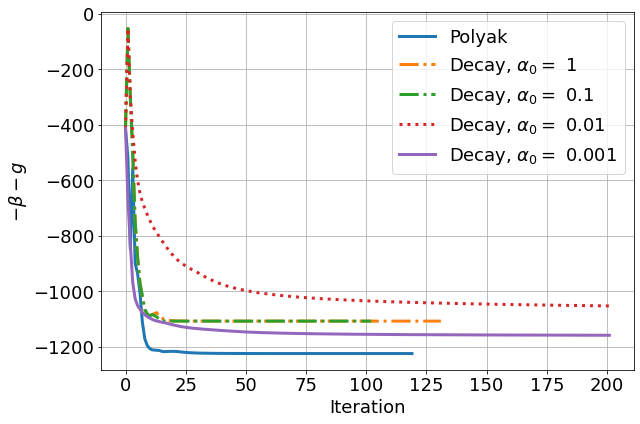}
			\label{fig:step-size-comp}
		\end{figure}	
	\end{minipage}
	\hspace{0.05\linewidth}
	\begin{minipage}{0.45\linewidth}
	\begin{figure}[H]
		\caption{Algorithm \ref{alg:pricing} converges to different local solutions for different starting points}
		\includegraphics[width=\linewidth]{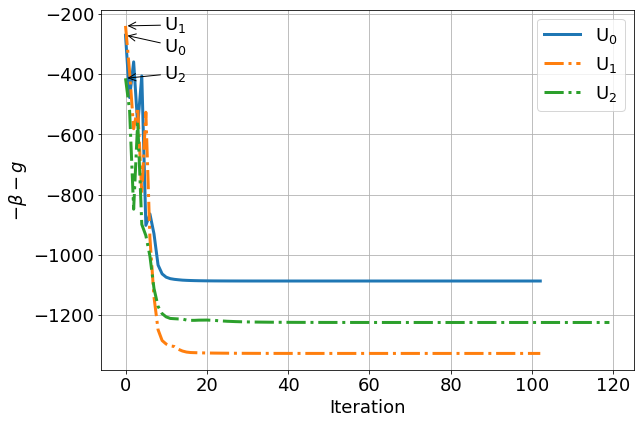}
		\label{fig:multi-start}
	\end{figure}
	\end{minipage}
\end{minipage}

\section{Comparison with Existing Methods}
\label{sec:experiments}
We next compare the performance of \ourmethod against various methods from the literature. 

 
\subsection{Synthetic Dataset}
The difficulty of SCMD depends on several factors such as the arrangement of subspaces, the separation between
subspaces, the total dimensions of the subspaces, and the percentage of the missing data. 
In this section, in addition to the random instances discussed in Section~\ref{subsec:experimental-setup}, we also
consider instances generated from disjoint\footnote{A collection of subspaces $S_1,\dots, S_c $ is said to
be disjoint if $\text{dim}(S_i)+\text{dim}(S_j)=\text{dim}(S_i\cup S_j)$ and $S_i\cap S_j =\{0\}\ \forall i\neq j$).} semi-random subspace arrangements. Semi-random instances allow us to control the separation between the subspaces which is measured by the affinity. 

We generate intances having two and three disjoint subspaces
with affinity between them being controlled by an angle parameter
$\theta \in [0,\frac{\pi}{2}]$.  Small values of $\theta$ indicate low affinity between the subspaces, and hence the
clustering task is more challenging \citep{subsace-clustering-geometric-outliers}. 
Similar to \citep{review-maryam-2021}, we generate the two disjoint subspaces as follows:
\[\basis[1]=\Big( {I_{\lowrank} \atop 0_{\lowrank}} \Big), 
\basis[2]=\Big( {\cos(\theta)I_{\lowrank} \atop \sin(\theta)I_{\lowrank}} \Big)\]
Here $\lowrank[1]=\lowrank[2]=\lowrank$, $I_{\lowrank}$ denotes the identity matrix of size  $\lowrank \times \lowrank$,
and $0_{\lowrank}$ denotes the zero matrix of size $\lowrank \times \lowrank$.  
For an instance with $\numvectors$ data points, we generate $\frac n2$ data points from each of the two subspaces. 
We first randomly create data points within each of the two $2\lowrank$-dimensional subspaces and then transform them to the $\ambientdimension$-dimensional space.
Let $\hat{X}_1$ and $\hat{X}_2$ denote data points generated from each subspace within $2r$-dimensional space. 
$\hat{X}_1$ is created from $\basis[1]$ as $\hat{X}_1=\basis[1]\times W$ where $W\in \mathbb{R}^{\lowrank\times
\frac{\numvectors}{2}}$ and each entry of $W \sim \mathcal{N}(0,1)$. 
$\hat{X}_2$ is created similarly from $\basis[2]$.  
These data points from dimension $2\lowrank$ are then transformed to ambient dimension $\ambientdimension$ by
multiplying with a randomly generated orthonormal basis $P\in R^{d\times2\lowrank}$ as $X_i=P\times \hat{X}_i$ for $i=1,2$. 
This orthonormal projection preserves the affinity between two subspaces from the $2\lowrank$-dimensional space
to the $\ambientdimension$-dimensional space. 

We also create instances with three disjoint subspaces, and generate $\numvectors/3$ vectors from each of the subspaces.  The three initial subspaces are constructed as \citep{review-maryam-2021}:
\[
\basis[1]=\Big( {I_{\lowrank} \atop 0_{\lowrank}} \Big), 
\basis[2]=\Big( {\cos(\theta)I_{\lowrank} \atop \sin(\theta)I_{\lowrank}} \Big),
\basis[3]=\Big( {-\cos(\theta)I_{\lowrank} \atop -\sin(\theta)I_{\lowrank}} \Big)
\]
The rest of the construction is identical to the two disjoint subspace case.

After generating data matrix $X$ either randomly or via the disjoint subspaces approach,  we uniformly at random drop a percentage $\missingpercentage$ of the entries of the matrix $X$ yielding the set of observed entries $\observedset[]$.
\subsection{Metrics}
We compare  performance of all the methods in terms of clustering error and completion error defined as follows:
\begin{itemize}
	\item Clustering Error: Let $\{G_1,G_2,\dots,G_\numsub\}$ be the ground truth clusters (which are known for our
	synthetic instances) where $ G_k \subseteq [\numvectors]\ \forall k\in[\numsub]$ and similarly let $\{P_1,P_2,\dots,P_{\numsub'}\}$ be the predicted clusters. 
	We evaluate the predicted clusters by solving the following assignment problem to find the best matching between predicted and true clusters: 
	\begin{subequations}
		\label{eq:opt-matching} 
		\begin{alignat}{2}
			\hat{c} := \min_{y\in \{0,1\}^{\numsub \times \numsub'}}  & \sum_{k \in [\numsub]} \sum_{k' \in [\numsub']}|P_{k'}\triangle
			G_k| y_{kk'} \label{eq:opt-matching-obj}\\
			&\sum_{k' \in [\numsub']} y_{kk'} = 1 && \forall k \in [\numsub] \label{eq:opt-matching-con1}\\
			&\sum_{k \in [\numsub]} y_{k k'} \leq 1 && \forall k' \in [\numsub'] \label{eq:opt-matching-con2}
		\end{alignat}
	\end{subequations}
where binary variable $y_{kk'}=1$ if and only if true cluster $k \in [\numsub]$ is matched with predicted cluster $k' \in
[\numsub']$ and the objective coefficients  $|P_{k'}\triangle G_k|$ measure the size of the symmetric difference between
predicted cluster $P_{k'}$ and true cluster $G_k$.
Thus, the objective \eqref{eq:opt-matching-obj} minimizes the number of disagreements in the matched clusters,
constraints~\eqref{eq:opt-matching-con1} ensure that each true cluster is mapped to
exactly one predicted cluster, while constraints \eqref{eq:opt-matching-con2} require that each predicted cluster is
mapped to at most one true cluster. 
Formulation \eqref{eq:opt-matching} is valid when $K' \geq K$ -- in case $K' < K$, the equations in
\eqref{eq:opt-matching-con1} are switched to $\leq$ and the inequalities \eqref{eq:opt-matching-con2} are switched to
equations. After solving \eqref{eq:opt-matching}, 
we can calculate the clustering error as  $100\times\hat{c}/2\numvectors$. 
	Note that we divide by 2 since every vector is penalized twice when mismatched. 
	
	\item Completion Error: Let $\observedset^c$ denotes the set of missing indices and $I_{\observedset^c}$ be the
	projection operator restricted to  $\observedset^c$. We define completion error to be the relative Frobenius distance
	between the true and recovered unobserved entries: $\|I_{\observedset^c}\circ(\hat{X}-X_{GT})\|_F/\|I_{\observedset^c}\circ(X_{GT})\|_F$.
	Here $\hat{X}$ refers to completed matrix and $X_{GT}$ refers to the ground truth matrix. Once we recover the
	clusters, we perform low-rank matrix completion on the data corresponding to each cluster separately to construct
	$\hat{X}$. This step is done using GROUSE \citep{laura-grouse} if the subspace dimensions are assumed known and singular value thresholding (SVT)
	\citep{lrmc-svd} if the dimensions are unknown.
\end{itemize}
\subsection{Comparison against other MIP approaches}
\label{subsec:comp-mip}
We first benchmark \ourmethod against the following MILP based facility-location methods proposed in the literature. 
These methods were proposed for the fully-observed data case. We do natural extensions to account for missing data as follows:
\begin{itemize}
	\item FLoSS \citep{floss}:  
	The candidate subspaces in FLoSS are initialized from the data by randomly selecting $\lowrank$ sets of linearly independent points, with $2\leq \lowrank < \ambientdimension$. 
	Data corresponding to  each set of $\lowrank$ points defines a linear subspace of dimension ($\lowrank-1$). 
	The corresponding basis for fully-observed data is obtained by performing SVD on the sampled points to get the best fit subspace $\basis$. 
	We extend this approach to partially observed data by using the same cost model as ours, i.e., the residual on observed entries~\eqref{eq:residual-closed-form}. 
	To handle missing data in the subspace generation process, we perform LRMC using GROUSE on each set of sampled points to get the best fit subspace.
	However, instead of sampling $\lowrank$ vectors, we sample $2\lowrank$ vectors  since LRMC is likely to fail with
	only $\lowrank$ vectors.  We then solve model~\eqref{eq:milp} with $[\numcols]$ consisting of all the subspaces generated with the above strategy. 
	We refer to this algorithm as MIP-RANDOM. 
	
	We point out that authors in~\citep{floss} do not specify the number of candidate subspaces ($|\numcols|$) to construct. 
	Since the candidate subspace generation process is random, MIP-RANDOM also serves as a good benchmark for \ourmethod
	to demonstrate the value of our approach for dynamically generating subspaces.
	Thus, we consider a high number of candidate subspaces, $|\numcols|=5000$, in MIP-RANDOM whereas we initialize \ourmethod with only 300 subspaces.
	
	\item MB-FLoSS \citep{minimal-basis-facility}: MB-FLoSS is similar to FloSS, with the difference being the candidate
	subspaces generation strategy. Instead of doing random sampling,  candidate subspaces are generated by solving the following optimization problem:
	\begin{alignat}{2}
		C^*=\argmin_C \|C\|_{2,1}+\gamma\|E\|_{1,2} \quad \text{s.t. } X=XC+E.
		\label{eq:mb-floss}
	\end{alignat}
	Here $C$ is the coefficient matrix of variables with $\|C\|_{2,1} := \sum_{i=1}^\ambientdimension\sqrt{\sum_{\datavector=1}^\numvectors (C_{i\datavector})^2}$, and $E$ is the error matrix of variables with $\|E\|_{1,2} := \sum_{\datavector=1}^\numvectors\sqrt{\sum_{i=1}^\ambientdimension (E_{i\datavector})^2}$.
	To extend this to the missing data case, we zero-fill the missing entries when solving~\eqref{eq:mb-floss} to get $C^*$.
	Each column $\datavector$ of $C^*$ represents the coefficients of other data points required to represent the data vector $j$. 
	With an estimate of subspace dimension $\lowrank$ at hand, data point $\datavector$  needs at most $\lowrank$ other
	points for representation. Thus, for each column, we use the data points associated with the $\lowrank$ largest
	coefficients in absolute value in that column to form a candidate subspace. 
	For fully-observed data, SVD is used on these data points to get a candidate subspace. 
	With missing data, we perform LRMC using GROUSE to get the best fit subspace which is then used as a candidate subspace.
	The number of candidate subspaces generated by this method is therefore the number of unique subspaces generated by
	considering each column of $C^*$. 
	
	\item BB-LRR \citep{hu-milp-2015}: 
	BB-LRR generates candidate subspaces by over-segmentation in low rank representation \citep{LRR}.
	They set the number of clusters larger than the ground truth in the spectral clustering step, specifically they use $\numsub+3$ instead of $\numsub$.
	For the assignment cost function, instead of using the the residual distance directly ($\distance$ as
	in~\eqref{eq:residual-closed-form}), they use a normalized distance that depends on a dimension-dependent goodness of fit of subspace.
	The authors also suggest a randomized local method for subspace generation but we don't benchmark against it for two reasons: a) this subspace generation process is very similar to FLoSS discussed above, and b) LRR generated subspaces outperformed randomized local models in \citep{hu-milp-2015}.
	We extend this algorithm to missing-data case by zero-filling the missing data during LRR step and doing LRMC with GROUSE for each cluster when generating a candidate subspace.

\end{itemize}
\paragraph{Disjoint subspaces:}
\begin{figure}
	\caption{Performance comparison for different MIP-based methods as a function of subspace angles for two disjoint subspaces . Parameters are $\mathbf{\ambientdimension=20, \numvectors =200},\mathbf{\numsub}=2$, and $\mathbf{\lowrank[1] =\lowrank[2]=2}$}	
	\label{fig:indep-ss-mip}	
	\centering
	\begin{subfigure}{0.85\linewidth}
		\centering
		\caption{Clustering error}
		\includegraphics[width=\linewidth]{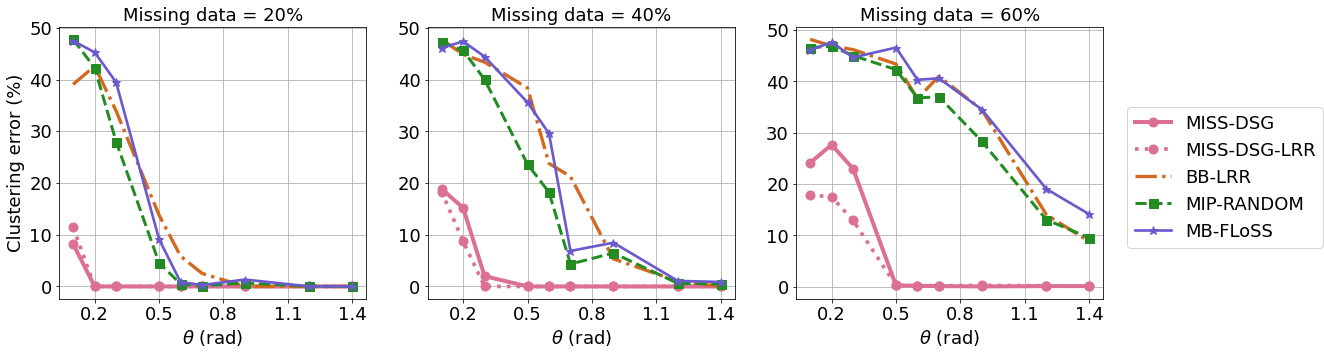}
		\label{fig:indep-ss-mip-clustering}	
	\end{subfigure}
	\begin{subfigure}{0.85\linewidth}
		\centering
		\caption{Completion error}
		\includegraphics[width=\linewidth]{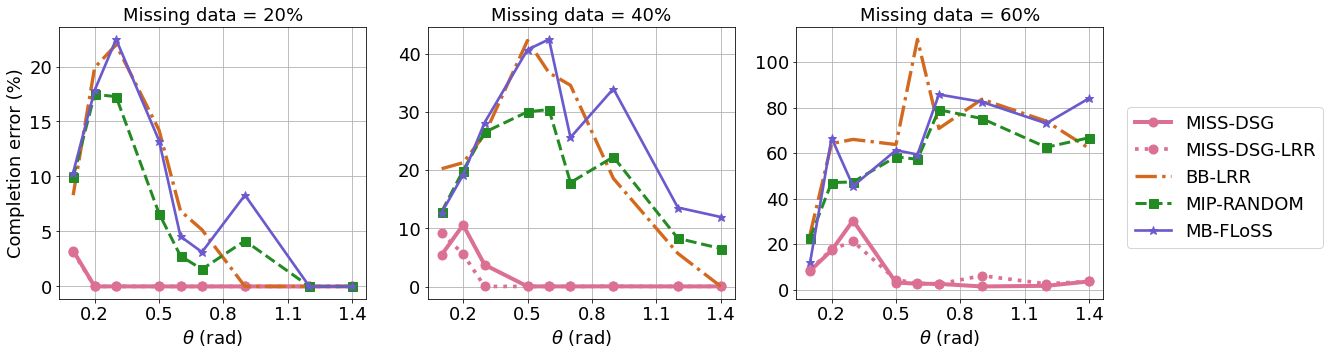}
		\label{fig:indep-ss-mip-completion}	
	\end{subfigure}
\end{figure}
We fix parameters $\ambientdimension=20$, $\numvectors =200$, and $\lowrank[i] =
2 \ \forall \ i = 1, 2$. We vary $\theta$ between $0.1 (\approx 6^\circ)$ and $1.4(\approx 80^\circ)$. For each value of $\theta$, we consider 10 random trials for each setting. We let the missing data percentage $\missingpercentage \in \{20,40,60\}$\%, and we report clustering and completion errors in Figure \ref{fig:indep-ss-mip}. 
For our method, we consider both a random initialization (referred to as \ourmethodns) and an initialization with LRR similar to \cite{hu-milp-2015} (referred to as \ourmethodns-LRR).
Existing facility location methods give similar performance for  all three missing data cases as shown in Figure \ref{fig:indep-ss-mip-clustering}. 
We observe that the perfect recovery threshold in terms of clustering error for existing MIP-based methods is $\theta=0.5$ for $\missingpercentage=20\%$, $\theta=1.2$ for $\missingpercentage=40\%$, and $\theta>1.4$ for $\missingpercentage=60\%$. 
Thus, existing MIP-based methods fail when subspaces are in close affinity or there is a high amount of missing data
while \ourmethod still gives low clustering errors in this regime.  
We observe a similar trend in completion error as shown in Figure~\ref{fig:indep-ss-mip-completion}.
For $\missingpercentage=60\%$, we observe that for the existing methods the completion errors increases with $\theta$ while one expects it to
decrease. The candidate subspace generation strategy in these methods use self-expressiveness and low rank
matrix completion, both of which fail with high levels of missing data.
The LRMC step for calculating completion error can often be faulty when an estimated cluster has points from multiple subspaces, translating to arbitrary recovery of missing entries and hence high completion errors.
We also point out that \ourmethod added between 800-4000 new candidate subspaces in these instances while MIP-RANDOM was
initialized with 5000 candidate subspaces. 
Thus, the superior performance of \ourmethod over MIP-RANDOM demonstrates the value of solving pricing
problem~\eqref{eq:pricing} to generate new candidate subspaces.
We also highlight that BB-LRR and \ourmethodns-LRR are initialized with the same set of initial subspaces. 
\ourmethodns-LRR outperforms BB-LRR by a significant margin, demonstrating that even with non-random initial subspaces,
our method of generating candidate subspaces by solving pricing problem~\eqref{eq:pricing} leads to significant
improvements. 

\begin{figure}[]
	\centering
	\caption{Performance comparison for different MIP-based methods on randomly sampled subspaces}
	\begin{subfigure}{0.49\linewidth}
		\caption{Effect of missing data on clustering error}
		\label{fig:missing-data-mip}
		\includegraphics[width=0.85\linewidth]{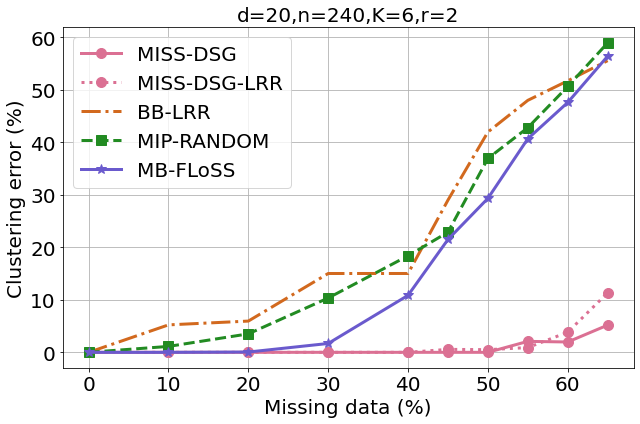}
	\end{subfigure}
	\begin{subfigure}{0.49\linewidth}
		\caption{Effect of $\ambientdimension/\numsub\lowrank$ on clustering error}
		\label{fig:rank-ratio-mip}
		\includegraphics[width=0.85\linewidth]{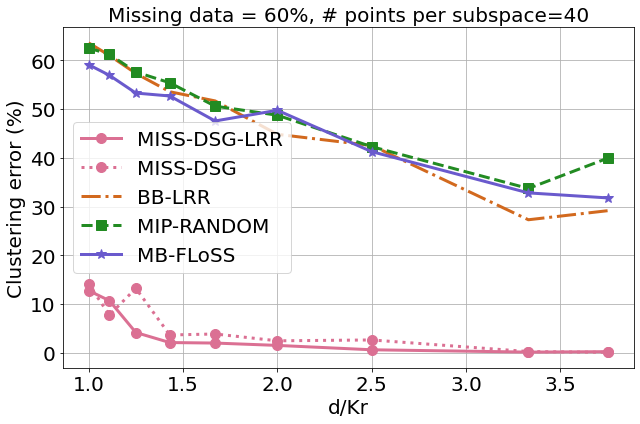}
	\end{subfigure}
\end{figure}

\paragraph{Random subspaces:}
We now consider randomly generated subspaces to study the effect of missing data and ratio of ambient dimension to total
rank ($\ambientdimension/\sum_{i=1}^\numsub \lowrank[i]$). 
We fix $\ambientdimension=20, \numvectors=240,\numsub=6,\lowrank[i]=r=2\ \forall i=1,\dots,K$, and vary the percentage
of missing data between 0 to 65\% as shown in Figure \ref{fig:missing-data-mip}. 
The reported results are averaged over 10 random trials.
We observe that \ourmethod yields low clustering error over a much wider range of missing data percentanges than the
other MILP methods. 
We report completion errors for these instances in Table \ref{tab:completion-missing-data}. \ourmethod recovers the
missing data with no completion error for $\missingpercentage$ up to $50\%$ while other MIP-based methods have high
completion errors for $\missingpercentage>20\%$.
We next study the effect of total rank of the data matrix on clustering error.
We generate a variety of instances with $\numvectors/\numsub\approx40, f=60\%$, and vary $\ambientdimension, \numsub, \lowrank$ to get $\ambientdimension/(\numsub\lowrank) \in[1,4]$. 
A lower ratio implies that the matrix is high-rank, thus making the clustering task more difficult. 
Due to a high amount of missing data, we observe high clustering errors in all the existing MIP-based methods. As the
matrix rank gets smaller relative to the ambient dimension, the clustering task becomes easier and hence the clustering error improves as shown in Figure \ref{fig:rank-ratio-mip}. 
\ourmethod recovers perfect clustering when the ratio is $>2$ and between $0-10\%$ for the ratio $\in [1,2]$. 
We observe that in the high-rank and high missing data regime, low clustering error does not imply low completion errors. 
This is  because low-rank matrix completion methods often fail in high-rank high missing data regime, leading to high completion errors, as shown in Table \ref{tab:completion-rank-ratio} for \ourmethod  for ratio $\in [1,1.7]$. 
\begin{table}[]
	\begin{minipage}{.48\linewidth}
		\caption{Average completion error  (\%) for random instances in Figure \ref{fig:missing-data-mip}}
		\label{tab:completion-missing-data}
		\centering
		\resizebox{\columnwidth}{!}{%
		\begin{tabular}{lrrrrr}
			\toprule
			$\missingpercentage$ & BB-LRR & MIP-RANDOM & MB-FLoSS & MISS-DSG & MISS-DSG-LRR \\
			\midrule
			10 & 28.3 & 2.9 & 0.0 & 0.0 & 0.0 \\
			20 & 43.6 & 14.5 & 0.6 & 0.0 & 0.0 \\
			30 & 116.9 & 43.7 & 17.1 & 0.0 & 0.0 \\
			40 & 126.2 & 104.3 & 80.7 & 0.0 & 0.0 \\
			50 & 237.8 & 198.1 & 173.0 & 0.1 & 6.7 \\
			55 & 267.4 & 259.7 & 257.7 & 35.2 & 12.7 \\
			60 & 324.0 & 311.4 & 295.5 & 41.9 & 43.0 \\
			65 & 368.7 & 361.3 & 338.5 & 114.5 & 154.4 \\
			\bottomrule
		\end{tabular}
	}
	\end{minipage}%
	\begin{minipage}{.5\linewidth}
		\centering
		\caption{Average completion error  (\%) for random instances in Figure \ref{fig:rank-ratio-mip}}
		\label{tab:completion-rank-ratio}
		\resizebox{\columnwidth}{!}{%
		\begin{tabular}{lrrrrr}
			\toprule
			$\ambientdimension/\numsub\lowrank$ & BB-LRR & MIP-RANDOM & MB-FLoSS & MISS-DSG& MISS-DSG-LRR \\
			\midrule
			1.1 & 428.3 & 437.8 & 399.9  & 102.1 & 142.5\\
			1.2 & 382.6 & 400.4 & 368.7 & 150.6 & 67.5  \\
			1.4 & 351.1 & 364.5 & 355.4  & 41.3 & 43.0\\
			1.7 & 324.0 & 311.4 & 295.5 & 43.0 & 41.9 \\
			2.0 & 286.4 & 297.7 & 296.1  & 42.4 & 64.7 \\
			2.5 & 223.4 & 226.4 & 225.8 & 33.1 & 31.6  \\
			3.3 & 140.5 & 187.3 & 165.4& 2.8  & 1.6 \\
			3.8 & 124.8 & 146.5 & 131.0 & 2.7 & 1.8  \\
			\bottomrule
		\end{tabular}
	}
	\end{minipage} 
\end{table}
\subsection{Comparison against state-of-the-art}
\label{subsec:comp-sota}
\begin{figure}
	\caption{Performance comparison against state-of-the-art methods on a variety of synthetic instances}
	\begin{subfigure}{0.47\linewidth}
		\caption{Effect of missing data on clustering error}
		\label{fig:missing-data-sota}
		\includegraphics[width=0.85\linewidth]{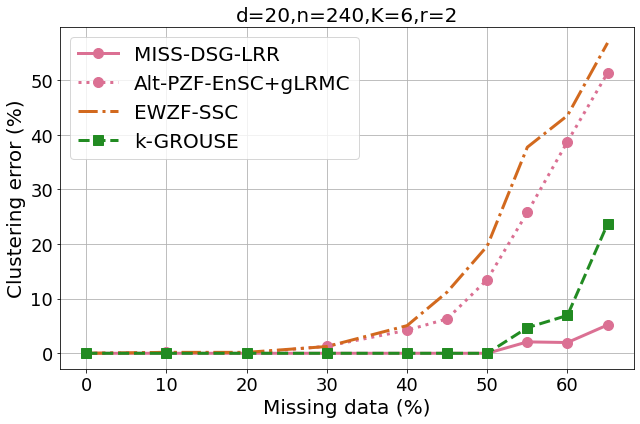}
		\end{subfigure}
	\begin{subfigure}{0.47\linewidth}
		\caption{Effect of $\ambientdimension/ \numsub \lowrank$ on clustering error}
		\label{fig:rank-ratio-sota}
		\includegraphics[width=0.85\linewidth]{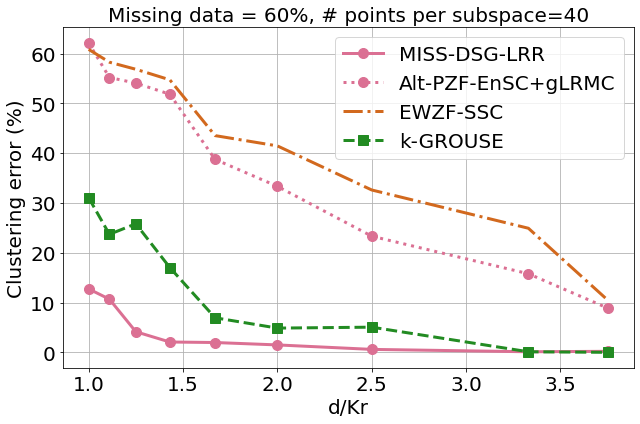}
	\end{subfigure}
\end{figure}

We now benchmark \ourmethod against the following methods from the literature:
\begin{itemize}
	\item EWZF-SSC: This is a natural extension of sparse subspace clustering to the case of missing data \citep{ssc_missing_icml2015}. In particular, \cite{ssc_missing_icml2015} proposed solving~\eqref{eq:ssc} with $\|\cdot\|_1$ regularization as follows:
	\begin{equation}
		C^*=\argmin \lambda\| I_{\observedset}\circ(X_{ZF}-X_{ZF}C)\|^2_F+ \|C\|_1 \qquad \text{s.t. diag}(C)=0
		\label{eq:ewzf-ssc}
	\end{equation}
	Coefficient matrix $C^*$ is then processed by the spectral clustering algorithm in order to obtain data segmentation as discussed in Section~\ref{sec:intro}.
	This method was found superior to the other methods proposed in \citep{ssc_missing_icml2015} for SCMD.
	\item \ssclrmcns: In a  review article on SCMD by \cite{Lane_2019_ICCV}, alternating between elastic-net subspace clustering \citep{ssc-ensc} and group low-rank matrix completion \citep{structured-ssc-vidal} was found to be the state-of-the-art method. PZF  is similar to EWZF and restricts error reduction on observed entries. The algorithm solves the following problem to get a coefficient matrix $C^*$, which is then processed by the spectral clustering algorithm to do data segmentation:
	\begin{equation}
		C^*=\argmin \lambda\| I_{\observedset}\circ(X_{ZF}-X_{ZF}C)\|^2_F+ \zeta \|C\|_1+(1-\zeta)\|C\|_F^2 \qquad \text{s.t. diag}(C)=0,
		\label{eq:ewzf-ensc}
	\end{equation}
	where $0<\zeta<1$.  The clusters obtained by spectral clustering algorithm are then processed group-wise by a low-rank matrix completion algorithm, e.g. SVT\citep{lrmc-svd}, to fill  the missing entries and get $\hat{X}$. In next iterations, $X_{ZF}$ in~\eqref{eq:ewzf-ensc} is replaced by $\hat{X}$, and algorithm alternates between clustering and completion for the given number of iterations.
	\item k-GROUSE: \cite{laura-k-grouse} proposed an extension of the K-Subspaces algorithm \citep{k-plane, q-flat} to the case of missing data. 
	The proposed algorithm is an alternating heuristic: starting with some initial subspaces, vectors are clustered by subspace assignment based on the same metric as~\eqref{eq:residual-closed-form}. 
	Given a cluster of vectors, matrix completion with Grassmannian Rank-One Subspace Estimation, GROUSE
	\citep{laura-grouse}, is performed to get a subspace estimate, and then vectors are reassigned, and the process is repeated until convergence. The algorithm stops when the clusters remain unchanged in successive iterations or the algorithm reaches the maximum allowed iterations
\end{itemize}
	Since \ourmethod and \ourmethodns-LRR gave similar performance, we report only
	\ourmethodns-LRR in these experiments.
	 The best parameter configurations are selected based on the average completion error on a hold out set. As noted by
	 \cite{Lane_2019_ICCV}, this approach translates more easily into practice compared to an  approach
	 where the parameter with least classification error is selected. Being an unsupervised learning task,  no true
	 cluster labels are available in practice, so the latter approach could not be implemented. However, one
	can hold out some observed entries as a validation set. In our experiments, we hold out 25\% of the data in
	the validation set for parameter selection.
	We report parameter choices for different algorithms in Table \ref{tab:parameter-choices}.
\begin{table}
	\caption{Parameter Choices}
	\label{tab:parameter-choices}
	\centering
	\begin{tabular}{ll}
		\hline
		\textbf{Method} & \multicolumn{1}{c}{\textbf{Parameter}}\\
		\hline \up
		\multirow{2}{*}{EWZF-SSC \citep{ssc_missing_icml2015}} & $\lambda = \frac{\alpha}{\max_{i\neq j} \|(X_{ZF})_{\observedset[j]}^T(X_{ZF})_{\observedset[j]})\|_{ij}}$\\ 
		&$\alpha \in \{5,20,50,100,200,320\} $\\
		k-GROUSE \citep{laura-k-grouse} & -\\
		\multirow{2}{*}{Alt-PZF-EnSC+gLRMC  \citep{Lane_2019_ICCV}}& $\lambda \text{ and } \alpha$ similar to EWZF-SSC\\\
		&$\zeta \in \{0.5,0.7,0.9\}$\\
		\hline
	\end{tabular}
\end{table}
\paragraph{Two disjoint subspaces:} 
\begin{figure}[]
	\centering
	\caption{Performance comparison against state-of-the-art as a function of subspace angles for two disjoint subspaces. Parameters are $\mathbf{\ambientdimension=20, \numvectors =200},\mathbf{\numsub=2}$, and $\mathbf{\lowrank[1] =\lowrank[2]=2}$
	}
	\label{fig:independent-ss-sota}
	\begin{subfigure}{0.85\linewidth}
		\centering
		\caption{Clustering error}
		\includegraphics[width=\linewidth]{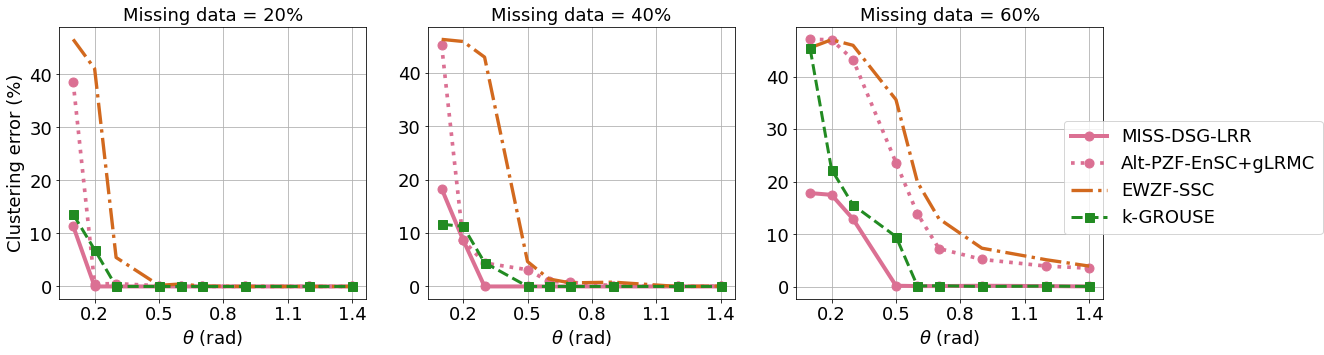}
		\label{fig:independent-ss-sota-clustering}
	\end{subfigure}
	\begin{subfigure}{0.85\linewidth}
		\centering
		\caption{Completion error}
		\includegraphics[width=\linewidth]{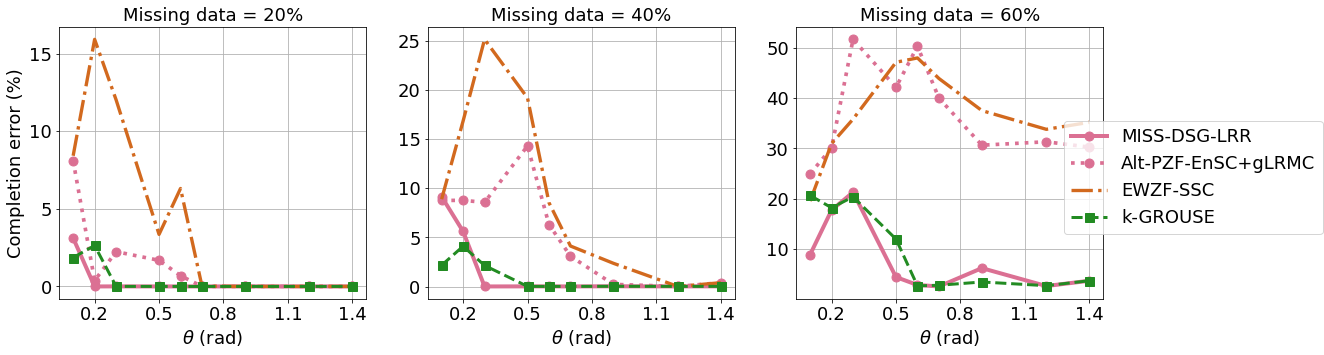}
		\label{fig:independent-ss-sota-completion}
	\end{subfigure}
\end{figure}
We keep the same experimental setting for disjoint subspaces as in Section \ref{subsec:comp-mip}. 
We report these results in Figure \ref{fig:independent-ss-sota}. 
As expected, EWZF-SSC fails when the percentage of the missing data is high or subspaces are close to each other (small $\theta$).
\ourmethodns-LRR and k-GROUSE give the lowest clustering errors. Both of these algorithms give similar performance with \ourmethodns-LRR doing slightly better than k-GROUSE for $\missingpercentage=60\%$.
A similar trend is observed in completion errors as shown in Figure~\ref{fig:independent-ss-sota-completion}.
\paragraph{Three disjoint subspaces:}  
\begin{figure}[]
	\caption{Performance comparison against state-of-the-art as a function of subspace angles for three disjoint subspaces. Parameters are $\mathbf{\ambientdimension=20, \numvectors =200},\mathbf{\numsub=3}$, and $\mathbf{\lowrank[1] =\lowrank[2]=\lowrank[3]=2}$
	}
	\label{fig:disjoint-ss-sota}
	\centering
	\begin{subfigure}{0.85\linewidth}
		\caption{Clustering error}
		\centering
		\includegraphics[width=\linewidth]{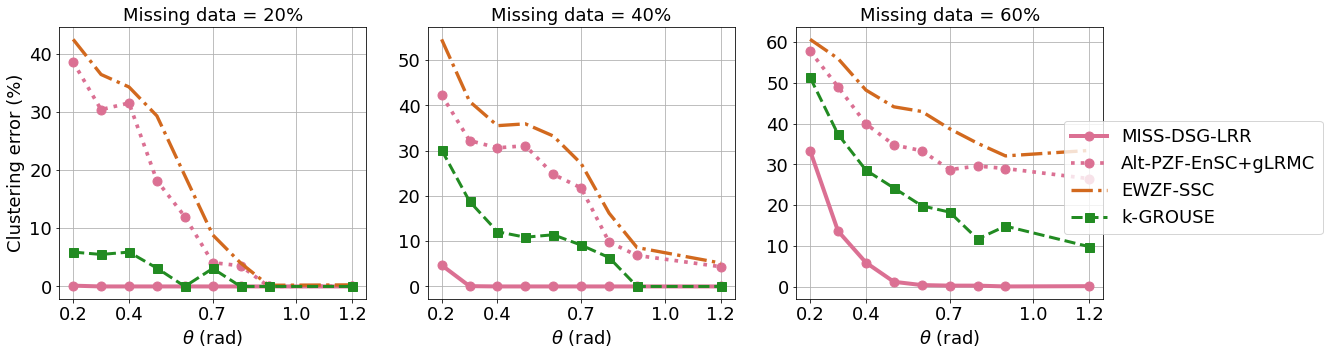}
		\label{fig:disjoint-ss-sota-clustering}
	\end{subfigure}
	\begin{subfigure}{0.85\linewidth}
		\centering
		\caption{Completion error}
		\includegraphics[width=\linewidth]{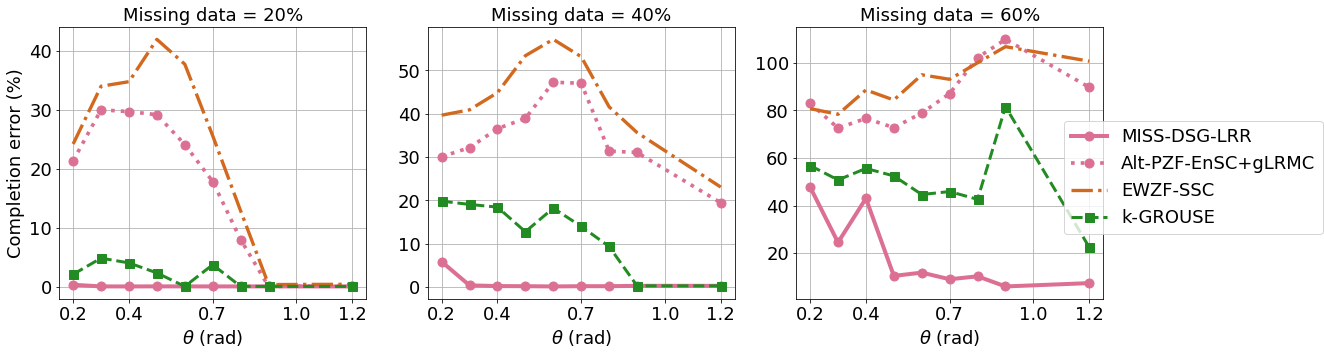}
		\label{fig:disjoint-ss-sota-completion}
	\end{subfigure}
\end{figure}
We fix parameters $\ambientdimension=20$, $\numvectors =200$, and $\lowrank[k] =
2 \ \forall \ k\in \{1,2,3\} $. We vary $\theta$ between $0.2\ (\approx 12^\circ)$ and $1.2\ (\approx68^\circ)$. For each value of $\theta$, we consider 10 random trials. We let missing data percentage $\missingpercentage \in \{20,40,60\}$. 
All methods are provided the true number of subspaces $\numsub$ and the true dimension of subspaces. 
We report these results in Figure \ref{fig:disjoint-ss-sota}. 
We see a significant drop in performance when compared to two disjoint subspaces for all algorithms except \ourmethodns-LRR.   
EWZF-SSC and \ssclrmc give high clustering errors when any pair of subspaces are close to each other (small $\theta$) or
the fraction of missing data is high, and are outperformed by both k-GROUSE and \ourmethodns-LRR.
Performance of k-GROUSE deteriorates in the low-affinity  and high missing data regimes.
\ourmethodns-LRR is the only algorithm which gives perfect recovery of clusters, in terms of low clustering errors as
well as low completion errors, in the low-affinity and high missing data regime.
\paragraph{Random subspaces:}
For random subspaces, we fix $\ambientdimension=20, \numvectors=240,\numsub=6,\lowrank[i]=r=2\ \forall i=1,\dots,K$, and
vary the percentage of missing data between 0 to 65\% as shown in Figure \ref{fig:missing-data-sota}. 
We observe that  EWZF-SSC and \ssclrmc exhibit significantly high clustering errors for $\missingpercentage>30\%$. 
k-GROUSE follows a similar trend with high clustering error for $\missingpercentage>50\%$.
In the high-missing data regime (40-65\%), \ourmethod-LRR yields the smallest clustering error.
Completion error follows a similar trend with both EWZF-SSC and \ssclrmc giving high completion error for $\missingpercentage>30\%$ as shown in Table \ref{tab:completion-missing-data-sota}. 
  \ourmethodns-LRR  gives lower completion error than k-GROUSE for $\missingpercentage>50\%$ while both give no completion error for $\missingpercentage<50\%$.

We next study the effect of total rank ($\numsub \lowrank$) relative to the ambient dimension $\ambientdimension$ on clustering error. We consider the same instances as we did in Section \ref{subsec:comp-mip} for randomly sampled subspaces and vary $\ambientdimension/\numsub \lowrank \in [1,4]$ as shown in Figure \ref{fig:rank-ratio-sota}. 
Since self-expressive methods do not perform well with high missing data (f = 60\%), we find that both EWZF-SSC and \ssclrmc give high clustering errors in all cases.
Performance of all algorithms improve as we move from the high-rank to low-rank regime.
In the high rank regime, ($1<\ambientdimension/\numsub \lowrank<2$), only \ourmethodns-LRR gives near perfect classification while k-GROUSE gives errors between 5-25\%.
Since we have high missing data, we observe that when the matrix is nearly full-rank ($\ambientdimension/\numsub
\lowrank <2$), all methods give high completion errors including  \ourmethodns-LRR, even though it has small clustering errors as shown in Table \ref{tab:completion-rank-ratio-sota}. This is due to the fact that low-rank matrix completion fails in such a setting. 

\begin{table}
	\begin{minipage}{.48\linewidth}
		\caption{Average completion error (\%) for random instances in Figure \ref{fig:missing-data-sota}. Column `Alt-'
		refers to method Alt-PZF-EnSC+gLRMC.}
		\label{tab:completion-missing-data-sota}
		\centering
		\resizebox{\columnwidth}{!}{%
			\begin{tabular}{lrrrr}
				\toprule
				$\missingpercentage$& EWZF-SSC & Alt- & k-GROUSE & MISS-DSG-LRR \\
				\midrule
				10 & 0.1 & 0.1 & 0.0 & 0.0 \\
				20 & 0.2 & 0.0 & 0.0 & 0.0 \\
				30 & 11.4 & 12.0 & 0.0 & 0.0 \\
				40 & 41.1 & 34.1 & 0.0 & 0.0 \\
				50 & 109.1 & 73.5 & 0.1 & 0.1 \\
				55 & 169.9 & 124.6 & 38.0 & 12.7 \\
				60 & 210.0 & 175.6 & 71.0 & 43.0 \\
				65 & 231.5 & 226.7 & 188.4 & 154.4 \\
				\bottomrule
			\end{tabular}
		}
	\end{minipage}%
	\begin{minipage}{.5\linewidth}
		\caption{Average completion error (\%) for random instances in Figure \ref{fig:rank-ratio-sota}. Column `Alt-'
		refers to method Alt-PZF-EnSC+gLRMC.}
		\centering
		\label{tab:completion-rank-ratio-sota}		
		\resizebox{\columnwidth}{!}{%
			\begin{tabular}{lrrrr}
				\toprule
				$\ambientdimension/\numsub\lowrank$& EWZF-SSC & Alt- & k-GROUSE & MISS-DSG-LRR \\
				\midrule
				1.1 & 275.6 & 235.5 & 225.6 & 142.5 \\
				1.2 & 253.1 & 234.8 & 238.2 & 67.5 \\
				1.4 & 238.1 & 200.6 & 143.4 & 43.0 \\
				1.7 & 210.0 & 175.6 & 71.0 & 41.9 \\
				2.0 & 174.1 & 137.8 & 46.2 & 64.7 \\
				2.5 & 143.0 & 117.7 & 45.1 & 31.6 \\
				3.3 & 106.4 & 90.9 & 1.4 & 1.6 \\
				3.8 & 60.1 & 52.2 & 0.0 & 1.8 \\
				\bottomrule
			\end{tabular}
		}
	\end{minipage} 
\end{table}

\subsection{Choice of penalty parameter in \ourmethod}
\label{subsec:penalty-choice}
\begin{figure}[H]
	\centering
	\caption{Effect of $\mathbf{\regularizer}$ (x axis, log scale) on clustering error}
	\label{fig:mip-convergence}
	\begin{subfigure}{0.3\linewidth}
		\caption{$\numsub$ known, $\lowrank$ unknown}
		\includegraphics[width=\linewidth]{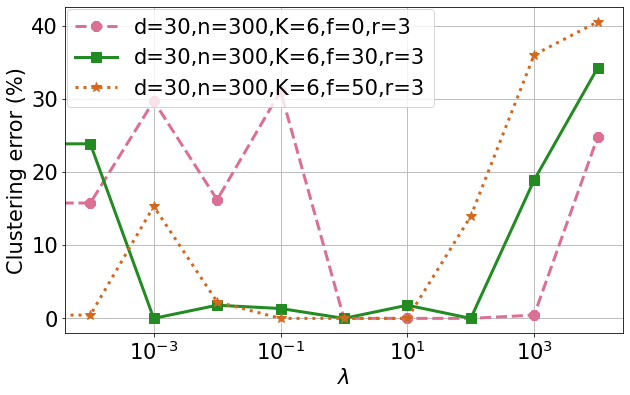}
	\label{fig:mip-convergence-K-1-r-0}
	\end{subfigure}
	\begin{subfigure}{0.3\linewidth}
		\caption{$\numsub$ unknown, $\lowrank$ known}
		\includegraphics[width=\linewidth]{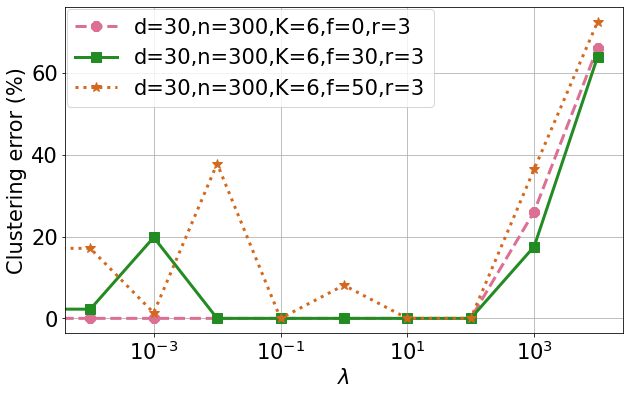}
	\label{fig:mip-convergence-K-0-r-1}
	\end{subfigure}
	\begin{subfigure}{0.3\linewidth}
		\caption{$\numsub$ unknown, $\lowrank$ unknown}
		\includegraphics[width=\linewidth]{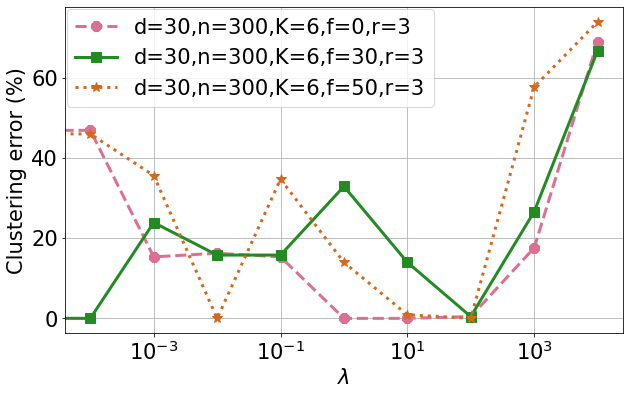}
	\label{fig:mip-convergence-K-0-r-0}
	\end{subfigure}
\end{figure}

If we know the number of subspaces ($\numsub$) and their underlying dimension $\lowrank[i]\ \forall i=1,\dots,\numsub$, then
we can use $\regularizer=0$ in the objective of~\eqref{eq:milp-ed}. 
If we do not have that information the choice of penalty parameter becomes an important hyperparameter. 
We conduct experiments investigating the impact of $\lambda$ on instances with $\ambientdimension=30,\numvectors=300,\numsub=6,\lowrank = \lowrank[k]=3\ \forall k \in [\numsub]$, and vary $\missingpercentage\in \{10,30,50\}$.
We consider three different cases for \ourmethodns:
\begin{itemize}
	\item $\numsub$ known, $\lowrank$ unknown: We use $\lowrank[\max]=2*\lowrank=6$ in Algorithm~\ref{alg:miss-dsg}. Thus our model considers subspaces of dimension $\in \{1,2,3,4,5,6\}$.
	However, we assume that we know number of subspaces ($\numsub$), and thus we keep constraint~\eqref{eq:beta-constr}.
	
	\item $\numsub$ unknown, $\lowrank$ known:  \ourmethod considers subspaces only of dimension $3$.
	 Since we don't know the number of subspaces ($\numsub$), we remove constraint~\eqref{eq:beta-constr} from our model, and let \ourmethod self-determine the number of subspaces.
	\item $\numsub$ unknown, $\lowrank$ unknown: 
	We again use $\lowrank[\max]=2*\lowrank=6$ in Algorithm~\ref{alg:miss-dsg}, and hence \ourmethod considers subspaces of dimension $\in \{1,2,3,4,5,6\}$.
	Similar to the previous case, we don't know $\numsub$ and hence constraint~\eqref{eq:beta-constr} is removed from the model.
	Thus, \ourmethod has freedom in selecting number of subspaces as well as their dimensions.
\end{itemize}
We report the effect of $\regularizer$  on clustering error in Figure \ref{fig:mip-convergence} for all three cases discussed above. 
We observe that \ourmethod gives low clustering errors for a wide range of $\regularizer$ values when either $\numsub$ or $\lowrank$ is known (Figures \ref{fig:mip-convergence-K-1-r-0}, \ref{fig:mip-convergence-K-0-r-1}). 
Extremely high values of $\regularizer$ lead to high clustering errors in all cases since the model is forced to either
select fewer subspaces or select subspaces of lower dimension than the ground truth. 
For $\regularizer=10^4$, \ourmethod selected subspaces of dimension $1$.
Similarly, for extremely small values of $\regularizer$, \ourmethod selects higher complexity subspaces, i.e., subspaces of dimension higher than ground truth if $\lowrank$ is not known or higher number of subspaces than ground truth if $\numsub$ is unknown.

Results for the case when we do not know either $\lowrank$ or $\numsub$ are shown in Figure~\ref{fig:mip-convergence-K-0-r-0}.
We observe that \ourmethod gives low clustering error only for a narrow range of $\regularizer$. 
This behavior is expected since the model has high degree of freedom. There are multiple union of subspace models which
might give  low assignment cost on the observed entries ($\sum_{\datavector\in[\numvectors]}\costvar$) but have a different complexity than ground truth. Hence, choice of $\regularizer$ is critical in this case.
\subsection{Hopkins155 Data Experiments}
Motion segmentation has been a standard dataset in the literature for benchmarking performance of SCMD algorithms.
Motion segmentation refers to the task of identifying multiple spatiotemporal regions corresponding to different rigid-body motions in a video sequence.
We consider Hopkins155 motion segmentation dataset \citep{Hopkins155} which contains 155 video sequences with 2 or 3 moving objects.
In each sequence, objects moving along different trajectories and all the trajectories associated with a single rigid
motion live in a 3-dimensional affine subspace \citep{ssc-ehsan-2009}. Similar to \cite{ssc_missing_icml2015}, we
subsample trajectories with six frames (equally spread) to simulate a high-rank data matrix. We handle affine subspaces
in our framework by considering an affine subspace of dimension $\lowrank$ in $R^{\ambientdimension}$ as a linear
subspace of dimension $\lowrank+1$ in $R^{\ambientdimension+1}$. Hence, we set $\lowrank[max]=4$ in our model. Since we
do not have information on the exact dimension of the underlying subspaces, we let our model self-determine it. However, for a fair comparison with other models, we do provide the number of subspaces as input to the model.  
We consider two variants of our methods: \ourmethod with random initialization, referred to as \ourmethod in Table
\ref{tab:hopkins} and \ourmethod with initialization from \ssclrmc, referred to as \ourmethodns-A in Table \ref{tab:hopkins}. For \ourmethodns-A, we let $\regularizer=0.1$ in our algorithm.
For state-of-the-art methods, we choose the best hyperparameter as discussed in Section \ref{subsec:comp-sota}. We report average clustering error over 155 sequences for each method for different missing data percentage in Table \ref{tab:hopkins}.
We observe that EWZF-SSC and k-GROUSE give similar performance with errors between 15-25\%  as missing data percentage is increased from 10\% to 50\%. \ourmethod gives error between 18-20\% for all values of missing data and is outperformed by \ssclrmc in all cases. However, we observed that initializing \ourmethod with \ssclrmc generated clusters offered a great advantage. With this initialization, \ourmethodns-A was able to improve upon \ssclrmc and gave errors between 5-13.9\%.
\begin{table}[h!]
	\caption{Clustering error (\%) by different methods on Hopkins 155 dataset with \# frames=6 and varying levels of
	missing data ($f$).}
	\label{tab:hopkins}
	\centering
	\resizebox{0.8\columnwidth}{!}{
		\begin{tabular}{lrrrrr}
			\toprule
			$\missingpercentage$& EWZF-SSC& k-GROUSE & \ssclrmc & \ourmethodns-A& \ourmethod \\
			\midrule
			10 &  17.4& 15.7&11.7& {5.3}&19.4\\
			20 &19.7 &  19.2  & 11.1     &  {5.8}&18.4\\
			30 & 20.1& 20.2      &  10.9    &    {6.5}& 20.2\\
			40 &20.2 &   21.5    &   12.6   &    {9.7}& 19.1\\
			50 &23.4 &25.1 &15.4  &  {13.9}&21.7  \\
			\bottomrule
		\end{tabular}
		}
\end{table}

\subsection{Computation Times}


We finally discuss the computational performance of \ourmethodns.
We first point out that most existing methods are faster than \ourmethodns\ -- the advantage of 
\ourmethodns\ is its ability to successfully cluster in cases where the other methods fail.
Hence, we only discuss computation times of existing methods briefly. 
We first discuss the computation times of existing MIP methods.
MB-FLoSS and BB-LRR are significantly faster than \ourmethod and MIP-RANDOM. 
Both BB-LRR and MB-FLoSS generate a small number of candidate subspaces using LRR and thus, solving the resulting MIP
model is typically fast ($<2$ minutes for synthetic instances considered in this paper). In fact, for these methods, the
majority of the time is spent in generating clusters using LRR. 
For MIP-RANDOM, we sample 5000 subspaces by performing LRMC on each sampled cluster of vectors. 
This step took between 4-25 minutes on average for the synthetic instances considered in this work. 
The subspace generation time varies based on the amount of missing data since LRMC also becomes expensive for instances
with high percentage of missing data. 
Due to a large number of candidate subspaces, solving the MIP model is also computationally more expensive than BB-LRR and MB-FLOSS, and took between 2-10 minutes.
EWZF-SSC was found to be computationally efficient on the considered synthetic instances, taking at most $2$ minutes. 
k-GROUSE took  a maximum of $12$ minutes.
\ssclrmc is more time-consuming with computational time varying between $3-120$ minutes. This includes the parameter training time considered in Table \ref{tab:parameter-choices}. For \ssclrmcns, we had a total of 18 choices for parameters tuning. 
For a single parameter choice, \ssclrmc gives similar computational efficiency as k-GROUSE.
\ourmethod required between $3-75$ minutes for the synthetic instances. 
We found that instances with a high-rank matrix required more computation time. 
Solving disjoint instances with \ourmethod was significantly faster than the random instances, taking a maximum of 5 minutes.

We refer readers interested in more details of the computation times for \ourmethod\ to the electronic companion to this paper.
\section{Conclusions and future directions}
We proposed a novel MILP framework \ourmethod for the Subspace Clustering with Missing Data problem and showed it is
capable of successfully clustering data in some regimes where all existing methods fail. \ourmethod offers several other potential advantages for SCMD. It gives the user flexibility to use a different function for cost of assignment between vector and subspace. If we know a good set of potential low dimensional subspaces, our framework can take advantage of this by including these subspaces in the formulation. \ourmethod is also capable of self-determining the number of subspaces and their dimensions, and can also easily be extended to include side constraints, e.g., ensuring that a given set of points does (or does not) lie in the same cluster. 
\ourmethod is computationally more expensive than the other clustering algorithms, and we leave speed improvements such as approximate gradient calculations in pricing and parallel implementation as future work.  

%
%
%
\ACKNOWLEDGMENT{Support for this research was provided by American Family Insurance through a research partnership with
the University of Wisconsin–Madison’s Data Science Institute.}
\def\urlprefix{}
\def\url#1{}
\bibliographystyle{informs2014} 
\bibliography{clustering-ijo} 

\begin{thebibliography}{46}
\providecommand{\natexlab}[1]{#1}
\providecommand{\url}[1]{\texttt{#1}}
\providecommand{\urlprefix}{URL }

\bibitem[{Abdolali \protect\BIBand{} Gillis(2021)}]{review-maryam-2021}
Abdolali M, Gillis N (2021) Beyond linear subspace clustering: A comparative
  study of nonlinear manifold clustering algorithms. \emph{Computer Science
  Review} 42:100435, ISSN 1574-0137,
  \urlprefix\url{http://dx.doi.org/https://doi.org/10.1016/j.cosrev.2021.100435}.

\bibitem[{{Balzano} et~al.(2010){Balzano}, {Nowak}, \protect\BIBand{}
  {Recht}}]{laura-grouse}
{Balzano} L, {Nowak} R, {Recht} B (2010) Online identification and tracking of
  subspaces from highly incomplete information. \emph{2010 48th Annual Allerton
  Conference on Communication, Control, and Computing (Allerton)}, 704--711,
  \urlprefix\url{http://dx.doi.org/10.1109/ALLERTON.2010.5706976}.

\bibitem[{{Balzano} et~al.(2012){Balzano}, {Szlam}, {Recht}, \protect\BIBand{}
  {Nowak}}]{laura-k-grouse}
{Balzano} L, {Szlam} A, {Recht} B, {Nowak} R (2012) K-subspaces with missing
  data. \emph{2012 IEEE Statistical Signal Processing Workshop (SSP)},
  612--615, \urlprefix\url{http://dx.doi.org/10.1109/SSP.2012.6319774}.

\bibitem[{Barnhart et~al.(1998)Barnhart, Johnson, Nemhauser, Savelsbergh,
  \protect\BIBand{} Vance}]{barnhart.et.al:98}
Barnhart C, Johnson EL, Nemhauser GL, Savelsbergh MWP, Vance PH (1998) Branch
  and price: Column generation for solving huge integer programs.
  \emph{Operations Research} 46:316--329.

\bibitem[{Benders(1962)}]{bnnobrs1962partitioning}
Benders JF (1962) Partitioning procedures for solving mixed-variables
  programming problems. \emph{Numerische mathematik} 4(1):238--252.

\bibitem[{Bradley \protect\BIBand{} Mangasarian(2000)}]{k-plane}
Bradley PS, Mangasarian OL (2000) k-plane clustering. \emph{Journal of Global
  Optimization} 16(1):23--32,
  \urlprefix\url{http://dx.doi.org/10.1023/A:1008324625522}.

\bibitem[{Cai et~al.(2010)Cai, Candès, \protect\BIBand{} Shen}]{lrmc-svd}
Cai JF, Candès E, Shen Z (2010) A singular value thresholding algorithm for
  matrix completion. \emph{SIAM Journal on Optimization} 20:1956--1982,
  \urlprefix\url{http://dx.doi.org/10.1137/080738970}.

\bibitem[{Candes \protect\BIBand{} Tao(2010)}]{lrmc-near-optimal}
Candes EJ, Tao T (2010) The power of convex relaxation: Near-optimal matrix
  completion. \emph{IEEE Transactions on Information Theory} 56(5):2053--2080,
  \urlprefix\url{http://dx.doi.org/10.1109/TIT.2010.2044061}.

\bibitem[{Candès \protect\BIBand{} Recht(2009)}]{recht_2009}
Candès EJ, Recht B (2009) Exact matrix completion via convex optimization.
  \emph{Foundations of Computational Mathematics} 9(6):717–772,
  \urlprefix\url{http://dx.doi.org/10.1007/s10208-009-9045-5}.

\bibitem[{{Charles} et~al.(2018){Charles}, {Jalali}, \protect\BIBand{}
  {Willett}}]{ssc-md-theroy-charles}
{Charles} Z, {Jalali} A, {Willett} R (2018) Sparse subspace clustering with
  missing and corrupted data. \emph{2018 IEEE Data Science Workshop (DSW)},
  180--184, \urlprefix\url{http://dx.doi.org/10.1109/DSW.2018.8439907}.

\bibitem[{Elhamifar(2016)}]{hrmc-elhamifar-NIPS2016}
Elhamifar E (2016) High-rank matrix completion and clustering under
  self-expressive models. Lee D, Sugiyama M, Luxburg U, Guyon I, Garnett R,
  eds., \emph{Advances in Neural Information Processing Systems}, volume~29
  (Curran Associates, Inc.),
  \urlprefix\url{https://proceedings.neurips.cc/paper/2016/file/9f61408e3afb633e50cdf1b20de6f466-Paper.pdf}.

\bibitem[{Elhamifar \protect\BIBand{} Vidal(2009)}]{ssc-ehsan-2009}
Elhamifar E, Vidal R (2009) Sparse subspace clustering. \emph{2009 IEEE
  Conference on Computer Vision and Pattern Recognition}, 2790--2797,
  \urlprefix\url{http://dx.doi.org/10.1109/CVPR.2009.5206547}.

\bibitem[{{Elhamifar} \protect\BIBand{} {Vidal}(2013)}]{ssc-Vidal-2013}
{Elhamifar} E, {Vidal} R (2013) Sparse subspace clustering: Algorithm, theory,
  and applications. \emph{IEEE Transactions on Pattern Analysis and Machine
  Intelligence} 35(11):2765--2781,
  \urlprefix\url{http://dx.doi.org/10.1109/TPAMI.2013.57}.

\bibitem[{Fan \protect\BIBand{} Chow(2017)}]{mc-fan-2017}
Fan J, Chow TW (2017) Matrix completion by least-square, low-rank, and sparse
  self-representations. \emph{Pattern Recognition} 71:290--305, ISSN 0031-3203,
  \urlprefix\url{http://dx.doi.org/https://doi.org/10.1016/j.patcog.2017.05.013}.

\bibitem[{Fischetti et~al.(2017)Fischetti, Ljubić, \protect\BIBand{}
  Sinnl}]{ufl-benders}
Fischetti M, Ljubić I, Sinnl M (2017) Redesigning benders decomposition for
  large-scale facility location. \emph{Management Science} 63(7):2146–2162,
  \urlprefix\url{http://dx.doi.org/10.1287/mnsc.2016.2461}.

\bibitem[{Ford \protect\BIBand{} Fulkerson(1958)}]{column-generation-ford}
Ford LR, Fulkerson DR (1958) A suggested computation for maximal
  multi-commodity network flows. \emph{Management Science} 5(1):97--101,
  \urlprefix\url{http://dx.doi.org/10.1287/mnsc.5.1.97}.

\bibitem[{{Hu} et~al.(2015){Hu}, {Feng}, \protect\BIBand{}
  {Zhou}}]{hu-milp-2015}
{Hu} H, {Feng} J, {Zhou} J (2015) Exploiting unsupervised and supervised
  constraints for subspace clustering. \emph{IEEE Transactions on Pattern
  Analysis and Machine Intelligence} 37(8):1542--1557,
  \urlprefix\url{http://dx.doi.org/10.1109/TPAMI.2014.2377740}.

\bibitem[{Huang et~al.(2004)Huang, Ma, \protect\BIBand{} Vidal}]{gpca-mix-UoS}
Huang K, Ma Y, Vidal R (2004) Minimum effective dimension for mixtures of
  subspaces: a robust gpca algorithm and its applications. \emph{Proceedings of
  the 2004 IEEE Computer Society Conference on Computer Vision and Pattern
  Recognition, 2004. CVPR 2004.}, volume~2, II--II,
  \urlprefix\url{http://dx.doi.org/10.1109/CVPR.2004.1315223}.

\bibitem[{Lane et~al.(2019)Lane, Boger, You, Tsakiris, Haeffele,
  \protect\BIBand{} Vidal}]{Lane_2019_ICCV}
Lane C, Boger R, You C, Tsakiris M, Haeffele B, Vidal R (2019) Classifying and
  comparing approaches to subspace clustering with missing data.
  \emph{Proceedings of the IEEE/CVF International Conference on Computer Vision
  (ICCV) Workshops}.

\bibitem[{Lazic et~al.(2009)Lazic, Givoni, Frey, \protect\BIBand{}
  Aarabi}]{floss}
Lazic N, Givoni I, Frey B, Aarabi P (2009) Floss: Facility location for
  subspace segmentation. \emph{2009 IEEE 12th International Conference on
  Computer Vision}, 825--832,
  \urlprefix\url{http://dx.doi.org/10.1109/ICCV.2009.5459302}.

\bibitem[{Lecun et~al.(1998)Lecun, Bottou, Bengio, \protect\BIBand{}
  Haffner}]{MNIST}
Lecun Y, Bottou L, Bengio Y, Haffner P (1998) Gradient-based learning applied
  to document recognition. \emph{Proceedings of the IEEE} 86(11):2278--2324,
  \urlprefix\url{http://dx.doi.org/10.1109/5.726791}.

\bibitem[{Lee \protect\BIBand{} Cheong(2013)}]{minimal-basis-facility}
Lee C, Cheong L (2013) Minimal basis facility location for subspace
  segmentation. \emph{2013 IEEE International Conference on Computer Vision
  (ICCV)}, 1585--1592 (Los Alamitos, CA, USA: IEEE Computer Society), ISSN
  1550-5499, \urlprefix\url{http://dx.doi.org/10.1109/ICCV.2013.200}.

\bibitem[{{Li} \protect\BIBand{} {Vidal}(2016)}]{structured-ssc-vidal}
{Li} C, {Vidal} R (2016) A structured sparse plus structured low-rank framework
  for subspace clustering and completion. \emph{IEEE Transactions on Signal
  Processing} 64(24):6557--6570,
  \urlprefix\url{http://dx.doi.org/10.1109/TSP.2016.2613070}.

\bibitem[{Liu et~al.(2010)Liu, Lin, \protect\BIBand{} Yu}]{LRR}
Liu G, Lin Z, Yu Y (2010) Robust subspace segmentation by low-rank
  representation. \emph{Proceedings of the 27th International Conference on
  International Conference on Machine Learning}, 663–670, ICML'10 (Madison,
  WI, USA: Omnipress), ISBN 9781605589077.

\bibitem[{Lu et~al.(2013)Lu, Feng, Lin, \protect\BIBand{}
  Yan}]{ssc-trace-lasso}
Lu C, Feng J, Lin Z, Yan S (2013) Correlation adaptive subspace segmentation by
  trace lasso. \emph{2013 IEEE International Conference on Computer Vision},
  1345--1352, \urlprefix\url{http://dx.doi.org/10.1109/ICCV.2013.170}.

\bibitem[{Lu et~al.(2012)Lu, Min, Zhao, Zhu, Huang, \protect\BIBand{}
  Yan}]{LSR}
Lu CY, Min H, Zhao ZQ, Zhu L, Huang DS, Yan S (2012) Robust and efficient
  subspace segmentation via least squares regression. Fitzgibbon A, Lazebnik S,
  Perona P, Sato Y, Schmid C, eds., \emph{Computer Vision -- ECCV 2012},
  347--360 (Berlin, Heidelberg: Springer Berlin Heidelberg), ISBN
  978-3-642-33786-4.

\bibitem[{Ng et~al.(2001)Ng, Jordan, \protect\BIBand{} Weiss}]{Ng01onspectral}
Ng AY, Jordan MI, Weiss Y (2001) On spectral clustering: Analysis and an
  algorithm. \emph{ADVANCES IN NEURAL INFORMATION PROCESSING SYSTEMS}, 849--856
  (MIT Press).

\bibitem[{Nguyen et~al.(2019)Nguyen, Kim, \protect\BIBand{} Shim}]{lrmc-review}
Nguyen LT, Kim J, Shim B (2019) Low-rank matrix completion: A contemporary
  survey. \emph{IEEE Access} 7:94215--94237,
  \urlprefix\url{http://dx.doi.org/10.1109/ACCESS.2019.2928130}.

\bibitem[{Panagakis \protect\BIBand{} Kotropoulos(2014)}]{ssc-ensc-music}
Panagakis Y, Kotropoulos C (2014) Elastic net subspace clustering applied to
  pop/rock music structure analysis. \emph{Pattern Recognition Letters}
  38:46--53, ISSN 0167-8655,
  \urlprefix\url{http://dx.doi.org/https://doi.org/10.1016/j.patrec.2013.10.021}.

\bibitem[{{Pimentel} et~al.(2014){Pimentel}, {Nowak}, \protect\BIBand{}
  {Balzano}}]{daniel-scmd}
{Pimentel} D, {Nowak} R, {Balzano} L (2014) On the sample complexity of
  subspace clustering with missing data. \emph{2014 IEEE Workshop on
  Statistical Signal Processing (SSP)}, 280--283,
  \urlprefix\url{http://dx.doi.org/10.1109/SSP.2014.6884630}.

\bibitem[{Pimentel-Alarc{\'o}n \protect\BIBand{}
  Nowak(2016)}]{PimentelAlarcn2016TheIR}
Pimentel-Alarc{\'o}n DL, Nowak R (2016) The information-theoretic requirements
  of subspace clustering with missing data. \emph{ICML}.

\bibitem[{{Pimentel-Alarcón} et~al.(2016){Pimentel-Alarcón}, {Balzano},
  {Marcia}, {Nowak}, \protect\BIBand{} {Willett}}]{daniels-gscc}
{Pimentel-Alarcón} D, {Balzano} L, {Marcia} R, {Nowak} R, {Willett} R (2016)
  Group-sparse subspace clustering with missing data. \emph{2016 IEEE
  Statistical Signal Processing Workshop (SSP)}, 1--5,
  \urlprefix\url{http://dx.doi.org/10.1109/SSP.2016.7551734}.

\bibitem[{Pimentel-Alarcón et~al.(2015)Pimentel-Alarcón, Boston,
  \protect\BIBand{} Nowak}]{lrmc-sampling-daniel}
Pimentel-Alarcón DL, Boston N, Nowak RD (2015) A characterization of
  deterministic sampling patterns for low-rank matrix completion. \emph{2015
  53rd Annual Allerton Conference on Communication, Control, and Computing
  (Allerton)}, 1075--1082,
  \urlprefix\url{http://dx.doi.org/10.1109/ALLERTON.2015.7447128}.

\bibitem[{Polyak(1987)}]{polyak}
Polyak B (1987) \emph{Introduction to optimization} (Optimization Software,
  Inc).

\bibitem[{Ramlatchan et~al.(2018)Ramlatchan, Yang, Liu, Li, Wang,
  \protect\BIBand{} Li}]{matrix-completion-recommendation}
Ramlatchan A, Yang M, Liu Q, Li M, Wang J, Li Y (2018) A survey of matrix
  completion methods for recommendation systems. \emph{Big Data Mining and
  Analytics} 1:308--323,
  \urlprefix\url{http://dx.doi.org/10.26599/BDMA.2018.9020008}.

\bibitem[{Rao et~al.(2010)Rao, Tron, Vidal, \protect\BIBand{}
  Yu}]{motion-segmentation-vidal-2010}
Rao S, Tron R, Vidal R, Yu L (2010) Motion segmentation in the presence of
  outlying, incomplete, or corrupted trajectories. \emph{IEEE transactions on
  pattern analysis and machine intelligence} 32:1832--45,
  \urlprefix\url{http://dx.doi.org/10.1109/TPAMI.2009.191}.

\bibitem[{Recht(2011)}]{lrmc-recht-bounds}
Recht B (2011) A simpler approach to matrix completion. \emph{J. Mach. Learn.
  Res.} 12(null):3413–3430, ISSN 1532-4435.

\bibitem[{Soltanolkotabi \protect\BIBand{}
  Candés(2012)}]{subsace-clustering-geometric-outliers}
Soltanolkotabi M, Candés EJ (2012) A geometric analysis of subspace clustering
  with outliers. \emph{The Annals of Statistics} 40(4):2195--2238, ISSN
  00905364, 21688966, \urlprefix\url{http://www.jstor.org/stable/41806533}.

\bibitem[{Tron \protect\BIBand{} Vidal(2007)}]{Hopkins155}
Tron R, Vidal R (2007) A benchmark for the comparison of 3-d motion
  segmentation algorithms. \emph{2007 IEEE Conference on Computer Vision and
  Pattern Recognition}, 1--8,
  \urlprefix\url{http://dx.doi.org/10.1109/CVPR.2007.382974}.

\bibitem[{Tsakiris \protect\BIBand{} Vidal(2018)}]{ssc-md-theory-Tsakiris}
Tsakiris M, Vidal R (2018) Theoretical analysis of sparse subspace clustering
  with missing entries. Dy J, Krause A, eds., \emph{Proceedings of the 35th
  International Conference on Machine Learning}, volume~80 of \emph{Proceedings
  of Machine Learning Research}, 4975--4984 (Stockholmsmässan, Stockholm
  Sweden: PMLR),
  \urlprefix\url{http://proceedings.mlr.press/v80/tsakiris18a.html}.

\bibitem[{Tseng(2000)}]{q-flat}
Tseng P (2000) Nearest q-flat to m points. \emph{Journal of Optimization Theory
  and Applications} 105(1):249--252.

\bibitem[{Wang et~al.(2019)Wang, Xu, \protect\BIBand{} Leng}]{ssc-lrr}
Wang YX, Xu H, Leng C (2019) Provable subspace clustering: When lrr meets ssc.
  \emph{IEEE Transactions on Information Theory} 65(9):5406--5432,
  \urlprefix\url{http://dx.doi.org/10.1109/TIT.2019.2915593}.

\bibitem[{Yang et~al.(2015)Yang, Robinson, \protect\BIBand{}
  Vidal}]{ssc_missing_icml2015}
Yang C, Robinson D, Vidal R (2015) Sparse subspace clustering with missing
  entries. Bach F, Blei D, eds., \emph{Proceedings of the 32nd International
  Conference on Machine Learning}, volume~37 of \emph{Proceedings of Machine
  Learning Research}, 2463--2472 (Lille, France: PMLR),
  \urlprefix\url{http://proceedings.mlr.press/v37/yangf15.html}.

\bibitem[{You et~al.(2016)You, Li, Robinson, \protect\BIBand{}
  Vidal}]{ssc-ensc}
You C, Li CG, Robinson DP, Vidal R (2016) Oracle based active set algorithm for
  scalable elastic net subspace clustering. \emph{2016 IEEE Conference on
  Computer Vision and Pattern Recognition (CVPR)}, 3928--3937,
  \urlprefix\url{http://dx.doi.org/10.1109/CVPR.2016.426}.

\bibitem[{Zapata et~al.(2007)Zapata, Gonzalez-Mora, la~Torre, Guil,
  \protect\BIBand{} Murthi}]{application-image-classification}
Zapata EL, Gonzalez-Mora J, la~Torre FD, Guil N, Murthi R (2007) Bilinear
  active appearance models. \emph{2007 11th IEEE International Conference on
  Computer Vision}, 1--8 (Los Alamitos, CA, USA: IEEE Computer Society),
  \urlprefix\url{http://dx.doi.org/10.1109/ICCV.2007.4409185}.

\bibitem[{Zhuang et~al.(2012)Zhuang, Gao, Lin, Ma, Zhang, \protect\BIBand{}
  Yu}]{ssc-lrr-zhuang}
Zhuang L, Gao H, Lin Z, Ma Y, Zhang X, Yu N (2012) Non-negative low rank and
  sparse graph for semi-supervised learning. \emph{2012 IEEE Conference on
  Computer Vision and Pattern Recognition}, 2328--2335,
  \urlprefix\url{http://dx.doi.org/10.1109/CVPR.2012.6247944}.

\end{thebibliography}


\addtocounter{table}{-1}
\refstepcounter{table}\label{LASTTABLE}
\end{document}